\documentstyle[12pt]{article}
\input xypic
\xyoption{all}
\LaTeXdiagrams
\newtheorem{teo}{Theorem}
\newtheorem{lem}{Lemma}
\newtheorem{prop}{Proposition}
\newtheorem{cor}{Corollary}
\newtheorem{defi}{Definition}

\title{Obstructions to Coherence: \\ Natural Noncoherent Associativity}
\author{Noson S. Yanofsky\thanks{e-mail: yanofsky@math.mcgill.ca}
\\ Department of Mathematics and Statistics\\ 
McGill University, Montreal, QC, Canada, H3A 2K6}
\begin{document}
\maketitle
\define\AC {{\bf {\cal AC}}}
\define\Abar{{\bf \bar{A}}}
\define\xbar{{\bf \bar{X}}}
\define\Xbar{{\bf \bar{\cal X}}}
\define\ahat{{\bf \hat{A}}}
\define\Ahat{{\bf \hat{\cal A}}}
\define\aChat{{\bf \hat{ A \cal C}}}
\define\Achat{{\bf \hat{\cal A}C}}
\define\chat{{\bf \hat{C}}}
\define\Chat{{\bf \hat{\cal C}}}
\define\xhat{{\bf \hat{X}}}
\define\Xhat{{\bf \hat{\cal X}}}
\define\an{{\bf A_n}}
\define\An{{\bf A_n}}
\define\ACk{{\bf {\cal AC}_k}}
\define\ACn{{\bf {\cal AC}_n}}
\define\cn{{\bf C_n}}
\define\Cn{{\bf \cal C_n}}
\define\xn{{\bf X_n}}
\define\Xn{{\bf \cal X_n}}
\define\tenAbar{\otimes_{\bf \bar{A}}}
\define\alphabar{ \bar{\alpha}}
\define\tenbar{ \bar{\otimes}}
\define\betabar{ \bar{\beta}}
\define\Bbar{{\bf \bar{B}}}
\define\Bprbar{{\bf \bar{B'}}}
\define\BB'bar{{\bf \bar{BB'}}}
\define\Bit{{\bf \overline{It(B)}}}
\define\Bitn{{\bf It(B)_n}}
\define\Bitpr{{\bf \bar{It(B')}}}
\define\Mbar{{\bf \bar{M}}}
\define\Topbar{{ \bar{Top}}}
\define\Botbar{{ \bar{Bot}}}
\define\tenBbar{\otimes_{\bf \bar{B}}}
\define\tenMbar{\otimes_{\bf \bar{M}}}
\define\betaBbar{\beta_{\bf \bar{B}}}
\define\alphaMbar{\alpha_{\bf \bar{M}}}
\define\Fsigbar{{\bf \bar{F \sigma}}}
\define\Fsign{{\bf (F \sigma)_n}}
\define\Hn{{\bf H_n}}
\define\Mn{{\bf M_n}}
\define\Tn{{\bf T_n}}
\define\homBnB{{\bf Hom_{{\tiny Cat}}(B}^n{\bf ,B)}}
\define\homnk{{\bf Hom_A(n,k)}}
\definemorphism{mapsto}\solid\tip\stop
\definemorphism{dashto}\dashed\tip\notip
\define \tenA{\otimes_{\bf A}}
\define \homAcat{{\bf Hom_{st \otimes}(A,Cat)}}
\define\lam{\lambda}
\define\OP{{\bf{\cal O}}}
\hrulefill
\begin{flushleft}\bf{Abstract} \\ \end{flushleft}
We study what happens when coherence
fails.
Categories with a tensor product and a natural associativity
isomorphism that does not necessarily satisfy the pentagon coherence
requirements (called associative categories) are considered.
Categorical versions
of associahedra where naturality squares commute and pentagons do not,
are constructed (called Catalan groupoids, $\An$).
These groupoids are used in the construction of the
free associative category. They are also used in the
construction of the
theory of associative categories (given as a 2-sketch). Generators and
relations are given for the fundamental group, $\pi(\An)$, 
of the Catalan
groupoids -- thought of as a simplicial complex. These groups are shown
to be more than just free groups. Each
associative category, $\bf B$, has related fundamental groups $\pi(\bf
B_n)$ and homomorphisms $\pi(P_n):\pi(\An) \longrightarrow \pi(\bf
B_n)$. If the images of the
$\pi(P_n)$ are trivial, i.e. there is only
one associativity path between any two objects, then the category is
coherent.
Otherwise the images of $\pi(P_n)$ are obstructions to coherence. Some
progress is made in classifying noncoherence of associative
categories.\\
\hrulefill
\section{Introduction} 
The history of coherence theory has its beginnings in homotopy theory.
In
1963, J. Stasheff \cite{Stasheff a} investigated the conditions in which
an H-space has a homotopy associative multiplication. At around the same
time, D.B.A. Epstein \cite{Epstein a} came across some associativity
questions while dealing with Steenrod operations. With these papers in
mind, S. Mac Lane then wrote his classic paper on coherence
\cite{Mac Lane a}. He abstracted the problem to the following
categorical
question: Given a category $\bf B$ and a tensor product on it $\otimes:
{\bf B \times B} \longrightarrow \bf B$ that is associative up to a
natural isomorphism $\beta_{A,B,C}: A \otimes (B \otimes C)
\longrightarrow (A \otimes B) \otimes C$ when is there a unique
canonical map between two specified formal combinations of objects? In
other words, what conditions on $\bf B, \otimes$ and $\beta$ insure that
all combinations of identities and $\beta$ between any two objects are
the same. Mac Lane answered the question by giving the following
condition:
$$ \begin{diagram}
 & A((BC)D) \rrto^{\beta_{A,B \otimes C , D}}
& &(A(BC))D \drto^{\beta_{A,B,C} \otimes Id_D}\\
       A(B(CD))\urto^{Id_A \otimes \beta_{B,C,D}}
\ddrrto_{\beta_{A,B,C \otimes D}}& & & &((AB)C)D\\
\\
  & &(AB)(CD).\uurrto_{\beta_{A \otimes B,C,D}
}&
&
\end{diagram}$$
If this diagram commutes for every four objects $A, B, C$ and $D$, then
there is only one canonical map between any two objects.
An associativity
isomorphism that satisfies the above condition is called coherent and
categories with a coherent associativity isomorphism (and a unitary
requirement) were later
called monoidal or tensor categories.

Other coherence questions arose. There were
cartesian closed categories, monoidal closed categories, distributive
categories etc. In recent times, coherence problems have arisen in many
areas of mathematics and mathematical physics. Quantum
groups, quantum field theories, linear logic, knot theory etc. are just
a few of the diverse areas that now deal with coherence problems. (See
the last section for more examples and some references to
current coherence
problems.) Basically, categories with structure and coherence conditions
are seen as a type of higher dimensional algebras and such algebras
are ubiquitous in the mathematics done today.

What seems to have been left out is what happens if a coherence
condition
fails. Is there any structure that can still be recovered? Is there a
hierarchy of coherence conditions? The prototypical example of a
situation in
which coherence fails is when the category is ${\bf R-Mod}$, the
category of left
$R$ modules for an arbitrary ring R, the $\otimes$ is the usual tensor
product and the $\beta$ is the {\em un}usual $$\beta(a \otimes (b
\otimes c)) = -1 ((a \otimes b) \otimes c).$$ Given this $\beta$ the
above pentagon does not commute. Since $\beta$ is an isomorphism, we can
express this noncommutativity by saying that starting from $A \otimes
(B \otimes (C \otimes D))$ and going clockwise around the pentagon, we
do not get the identity map. However there is some structure left,
namely, going around the pentagon twice does give the identity. The
higher complexes do not commute. However, we may be able to say
something about the structure of the higher complexes.
The goal of this paper is to explore this structure.

We turn back to homotopy theory to study this higher structure. For
every natural number $n$, we construct a category $\An$ whose objects
correspond to
associations of $n$ letters and whose morphisms correspond to
reassociations. Since the reassociations are isomorphisms,
the $\An$'s are in fact groupoids. We call these groupoids the Catalan
groupoids. As with all groupoids, they can be thought of as simplicial
complexes: the objects are 0-cells, the morphisms are 1-cells and the
commuting parts will be 2-cells. One of the goals of this paper is to
calculate the fundamental groups of the $\An$'s and their quotients. If
a quotient of the $\An$'s is indiscrete (i.e. one morphism between any
two objects), then it is called coherent. If there is more then one
morphism, then we have an obstruction to coherence.

We begin by defining an associative category i.e.
a category with a bifunctor that is associative up to
a natural isomorphism that does not necessarily satisfy the
pentagon condition. Examples of such categories are given.
The Catalan groupoids, ( the $\An$s ) are constructed in section 3.
Section 4 uses the $\An$s to construct the free associative
category on one
generator, $\Abar$. The universal properties of $\Abar$ are proved.
Section 5 goes on
to show that all the $\An$s together have the structure of an
operad. This operad is used in the construction of $\bf A$, the 2-sketch
of the theory of associative categories. A short
discussion of associative categories that have units is given in section
6.

In section 8, the fundamental groups of each of the $\An$ is calculated.
To each $\An$ we assign a maximal indiscrete subgroupoid $\Tn$,
the categorical analogue of assigning a maximal tree to a simplicial
complex. The fundamental groups are then obtained by collapsing the 
$\Tn$ to a point . A way of describing the morphisms of $\An$ which are the
generators of the group is introduced. The general scheme for the
generators and relations are provided and the first seven groups are
presented. The seventh group is shown not to be a free group. All higher
groups are non-free groups. What follows is a discussion of
quotients of the 2-sketch of associative categories. An attempt is made
to classify the failure of coherence for both unital and nonunital
associative categories.

This paper ends with a section that lists some of the possible
applications of this work and ways we can go further in the study of the
failure of coherence.

I am grateful to my advisor, Prof. Alex Heller for many helpful
ideas and long
discussions. I would also like to thank my collegue M. Mannucci for
many stimulating conversations.
\section{Associative Categories}
\begin{defi}[Associative Category]
An Associative category is a category $\bf B$, a bifunctor
$\otimes_{\bf B} : \bf B \times \bf B \longrightarrow \bf B $
called ``tensor'', and a natural isomorphism
$$ \beta_{\bf B , \otimes} : \otimes_{\bf B} \circ ( Id_{\bf B} \times
\otimes_{\bf B} ) \longrightarrow \otimes_{\bf B} \circ (
\otimes_{\bf B} \times Id_{\bf B}) $$ i.e. for every $A,B,C$ in $\bf B$
an isomorphism $$\beta_{\bf B , \otimes,A,B,C }: A \otimes_{\bf B} (B
\otimes_{\bf B}
C ) \longrightarrow  (A \otimes_{\bf B} B ) \otimes_{\bf B} C$$ called
the ``reassociation''.
\end{defi}

We reserve the right to abandon the subscripts when there is no concern
for ambiguity. Discussions of a unit of the tensor will be left for
section 6. The important point is that we do not make any
coherence requirements.
\par
In general, all categories have a composition of morphisms
that is associative, however,
the name ``associative category'' is used here because the most
interesting
feature about our categories, is that their tensors are associative
up to a isomorphism.
\par

Examples of associative categories abound. Any monoidal category of
\cite{Mac Lane b} (also called a tensor category in the literature e.g.
\cite{Joy&Str})
is automatically an associative category. A noncoherent
example
of an associative category is ${\bf R-Mod}$, the category of $R$ modules
for a
commutative ring $R$. The tensor product is the usual tensor product of
modules and the reassociation
$$\beta_{A,B,C} : A \otimes (B \otimes C) \longrightarrow (A \otimes B)
\otimes C $$
is defined as
$$ a \otimes (b \otimes c) \longmapsto \zeta (a \otimes b) \otimes c$$
where $\zeta = -1$ (see \cite{Mac Lane a} ) or $\zeta$ is, for
example, the fifth root of unity. We shall come back to this example again
at the end of section 10.
\par
The above example can be abstracted and put into the language of
quantum groups. Let $A=(A,\Delta ,\varepsilon)$ be an algebra with a
comultiplication and a counit. Let $\Phi$ be an invertible element in
$A \otimes A \otimes A$ such that
$$ (Id \otimes \Delta)(\Delta(a)) = \Phi((\Delta \otimes Id)(\Delta(a)))
\Phi^{-1} ,$$
for all $a \in A$.
Then the category, ${\bf A-Mod}$, of $A$ modules, has the structure of an
associative
category. The tensor product of two modules is constructed using the
comultiplication and the associativity isomorphism is given as
$$ \beta_{A,B,C}(a \otimes (b \otimes c)) = \Phi ((a
\otimes b) \otimes c).$$
If $\Phi$ further satisfies
$$[(Id \otimes Id \otimes \Delta)(\Phi)][(\Delta \otimes Id \otimes
Id)(\Phi)] = [\Phi_{234}][(Id \otimes \Delta \otimes
Id)(\Phi)][\Phi_{123}]$$
where $\Phi_{123}=\Phi \otimes 1$ and $\Phi_{234}= 1 \otimes \Phi$, then
${\bf A-Mod}$ is, in fact, a coherent monoidal category (with proper
concern given to units). Such an algebra is called a Drinfeld algebra
\cite{Shn&Stern} or a quasi-bialgebra \cite{Kassel} (again, care must
be given to units.)

We will construct $\Abar$, the free associative
category on one generator. Roughly speaking, the objects in the
category will
be associations of n letters for any positive integer n (the elements of
the free ``anomic'' algebra with one binary operation -- sometimes
called a
``magma''
-- on one generator). Morphisms are
called reassociations. The tensor product, $\tenAbar$,
will concatenate two associations (ie. multiplication in the anomic
algebra). The tensor product of reassociations will be defined
similarly.

\section{$\An$, The Catalan Groupoids}
For each positive integer $n$, we will construct the
groupoid $ \An $ which has as objects associations
of $n$ letters and as morphisms reassociations.
The free associative category will
then be the disjoint union of all the $ \An$ i.e.
$$ \Abar = \coprod\limits_{n \in N^+}  \An . $$
These groupoids will be
called the Catalan groupoids. The Catalan
numbers
$$c_n =\frac{\left( \begin{array}{c} 2n-2 \\ n-1 \end{array} \right)
}{n}$$
are the number of associations of $n$ letters with no ambiguity in the
multiplication (see e.g. \cite{Hil&Ped}). These groupoids are the
categorical version of
what people who study finite complexes call associahedra (e.g.
\cite{Ziegler}).

The categories $\An$ are built up inductively in a manner not
unrelated to the
way Stasheff's complexes $K_n$ are built up (
\cite{Stasheff a}).
We let ${\bf A_1} = \bf 1 $, the trivial category with one object and
one identity morphism. Now assume that each of the $ {\bf A_1} $,$ {\bf
A_2} $, $ \ldots$, ${\bf A_{n-1}} $ is defined. We define ${\bf A_n}$
with the following pushout:
\begin{equation} \label{mainpo}
\begin{diagram}
\coprod\limits_{i+j+k = n} {\bf A_i \times A_j \times A_k \times
\dot{I}}
\ddto|<<\ahook
\rrto^{W_{i,j,k}} & &
\coprod\limits_{a+b = n} {\bf A_a \times A_b} \ddto^{U_{a,b}} \\
\\
\coprod\limits_{i+j+k = n} {\bf A_i \times A_j \times A_k \times I}
\rrto^{V_{i,j,k}} & &{\bf A_n}
\end{diagram}
\end{equation}
where $i,j,k,a$ and $b$ range over all positive integers and $\bf I$
(respectively
${\bf
\dot{I}}$) is the indiscrete
(resp.
discrete) category with two objects, $0$ and $1$. The left hand
vertical map is the obvious inclusion.
The map $W_{i,j,k}$ is defined for each $i,j,$ and $k$ as follows:
let $f,g,h$ be objects in ${\bf A_i}, {\bf A_j}, {\bf A_k}$ respectively
then $$W_{i,j,k}(f,g,h,t)= \left\{  \label{defofW}
\begin{array}{r@{\quad:\quad}l}
f, U_{j,k}(g,h) & t=0 \\
U_{i,j}(f,g), h & t=1
\end{array} \right. $$
where the $U_{j,k}$ and $U_{i,j}$ are defined from the previous
pushouts assumed by the induction hypothesis. $W_{i,j,k}$ is
defined for morphisms similarly.

A discussion of this pushout is in order.
Each association of n letters
must
have an ``outermost'' multiplication i.e. the last two segments of the
word to be associated. There are $a$ letters to the left of this
multiplication and $b$ letters to the right of this multiplication. Each
of these smaller words are also associated. Hence the category ${\bf A_a \times A_b}$.
This outermost multiplication can occur anywhere within the word,
hence the coproduct. All these smaller categories also have
reassociations and they are carried over to the to the new $\An$. There
are, however, new reassociations that are handled by the left-hand
side
of the pushout. Reassociations are concerned with three smaller words,
hence the ${\bf A_i \times A_j \times A_k}$ . There are two ways of
associating these three words into one whole word. $W_{i,j,k}$
maps these three words into the two ways of associating them. The
pushout connects these two ways with isomorphisms.

Let's look at the first few $\An$.
The objects of $\An$ are to be thought of as associations of n
letters. We write the letters of the associations with non-boldface
letters
$A,B,C
\ldots$ etc.
The letters should not be thought of as objects in a category. Rather
they are variables or ``place holders''. We write $A \otimes (B \otimes
C)$ as a shorthand for the functor $(-) \otimes ((-) \otimes (-)) =
\otimes \circ (Id \times \otimes)$ from a category cubed to itself.
             Reassociations will be written as $\alpha$. We also
show how each $\An$ is ``built up'' from the lower groupoids. This is
done by writing ${\bf A_a \times A_b}$ for its image under $U_{a,b}$.
Similarly for the image of ${\bf A_i \times A_j \times A_k \times I}$
under
$V_{i,j,k}$.

\begin{itemize}
\item ${\bf A_1}$ was defined to be $\bf 1$.
The reader can think of the category as looking like
$$ {\bf A_1} = A $$
where A is just a variable that corresponds to the single identity
functor.
\item ${\bf A_2}$ is defined from the following pushout:
$$
\begin{diagram}
{\bf \emptyset}
\ddto|<<\ahook
\rrto & &
{\bf 1 \times 1} = \bf 1 \ddto \\
\\
{\bf \emptyset}
\rrto & &{\bf A_2}
\end{diagram}
$$
The left side of the pushout is a vacuous coproduct. So
$$ {\bf A_2} = AB $$
\item  ${\bf A_3}$
$$ \begin{diagram}
\Text{ $ {\bf A_1 \times A_2} $ \\ $ A(BC) $ }
\xto[0,3]^{\bf A_1 \times A_1 \times A_1 \times I}_{\alpha}
&  &  & \Text{${\bf A_2 \times A_1} $\\$ (AB)C $}
\end{diagram}$$ 
\item ${\bf A_4}$
$$  \begin{diagram}
                    & A((BC)D) \rrto^{\bf A_1 \times A_2 \times A_1
 \times I}&       &(A(BC))D \drto^{\bf A_3 \times A_1}\\
        A(B(CD))\urto^{\bf A_1 \times A_3}
 \ddrrto_{\bf A_1 \times A_1 \times A_2 \times I}& & & &((AB)C)D\\
 \\
   & &\Text{${\bf A_2 \times A_2}$\\$(AB)(CD)$} \uurrto_{\bf A_2 \times A_1 \times A_1
 \times I}&
 &\\
 \end{diagram} $$
Note: This diagram does {\em not} commute.
\item  ${\bf A_5}$. Some of the names of the edges are
left out in order to make it more readable.
{\setlength{\oddsidemargin}{0mm}
\tiny
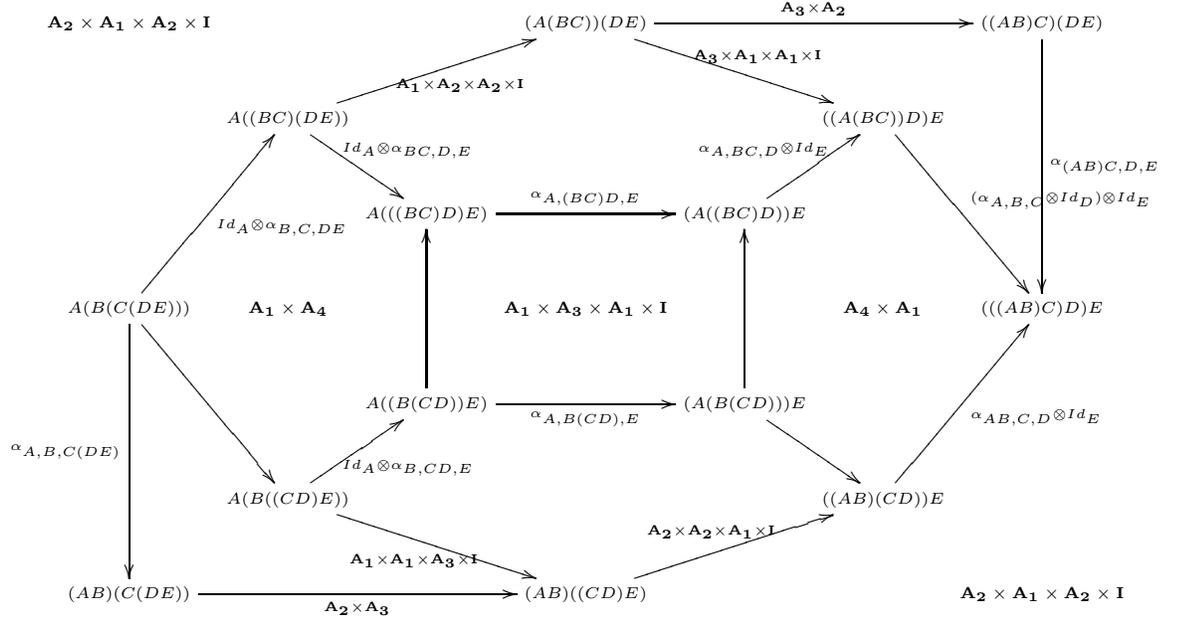
$$
\xymatrix@C=0pt{
{\bf A_2 \times A_1 \times A_2 \times I} &  &  &(A(BC))(DE)
\xto[0,3]^{\bf A_3 \times A_2}
\drrto^{\bf A_3 \times A_1 \times A_1 \times I}
&  &  &((AB)C)(DE)
\xto[3,0]^{\alpha_{(AB)C,D,E}}  \\
 &A((BC)(DE)) \urrto_{\bf A_1 \times A_2 \times A_2 \times I} \drto^{Id_A \otimes
\alpha_{BC,D,E}}
& & & &((A(BC))D)E
\ddrto^{(\alpha_{A,B,C} \otimes Id_D) \otimes Id_E}\\
 &  &A(((BC)D)E) \rrto^{\alpha_{A,(BC)D,E}} &  &(A((BC)D))E
\urto^{\alpha_{A,BC,D} \otimes Id_E}
& & \\
A(B(C(DE))) \uurto_{Id_{A}\otimes \alpha_{B,C,DE}}
\xto[3,0]_{\alpha_{A,B,C(DE)}}  \ddrto &{\bf A_1 \times A_4} & &{\bf A_1 \times A_3 \times A_1 \times I}  & &
{\bf A_4 \times A_1} &(((AB)C)D)E  \\
 & &A((B(CD))E) \rrto_{\alpha_{A,B(CD),E}} \uuto
& &(A(B(CD)))E
\uuto \drto  \\
 &A(B((CD)E)) \urto_{Id_{A}\otimes \alpha_{B,CD,E}}
\drrto_{\bf A_1 \times A_1 \times A_3 \times I }   &  &  &
&((AB)(CD))E
\uurto_{\alpha_{AB,C,D}\otimes Id_E} \\
 (AB)(C(DE)) \xto[0,3]_{\bf A_2 \times A_3}
& & &(AB)((CD)E)
\urrto^{\bf A_2 \times A_2 \times A_1 \times I}
& & & {\bf A_2 \times A_1 \times A_2 \times I}
}
$$}
\normalsize

${\bf A_1 \times A_1 \times A_3 \times I}$ in the lower left corner
describes the lower left quadrilateral. Similarly for ${\bf A_3 \times
A_1 \times A_1 \times I}$.
Notice the map $\alpha_{AB,C,DE}$ is shown twice: once going around the
top and once going around the bottom. There is only one such map and if
the page is bent, one will see that ${\bf A_5}$
is really a sphere.
\end{itemize}

\section{$\Abar$, The free associative category}
We now define the free associative category on one
generator.
$$ ( \Abar = \coprod\limits_{n \in N^+}  \An , \tenAbar,\alphabar).$$
The tensor product, $\tenAbar$, is defined as follows:
given $\phi :f \longrightarrow f'$ in ${\bf A_i}$,
$\gamma :g \longrightarrow g'$ in ${\bf A_j}$
and $\eta : h \longrightarrow h'$ in ${\bf A_k}$
then $f \tenAbar g = U_{i,j}(f,g)$ in ${\bf A_{i+j}}$.
The tensor of morphisms is defined as follows:
$\phi \tenAbar \gamma = U_{i,j}(\phi,\gamma)$ in ${\bf A_{i+j}}$.
Since $U_{i,j}$ is a functor, $\tenAbar$ is defined.
The reassociation is given by
$$\alphabar_{f,g,h} = V_{i,j,k}(Id_f,Id_g,Id_h, \iota):
f \tenAbar (g \tenAbar h) \longrightarrow (f \tenAbar g) \tenAbar h $$
where $\iota$ is the unique nontrivial isomorphism in $\bf I$.
Naturality means that the following diagram commutes
$$ \begin{diagram}
f \tenAbar (g \tenAbar h) \rrto^{\alphabar_{f,g,h}} \ddto_
{\phi \tenAbar( \gamma \tenAbar \eta )}& &
(f \tenAbar g) \tenAbar h \ddto^{(\phi \tenAbar \gamma ) \tenAbar \eta }\\
\\
f' \tenAbar (g' \tenAbar h') \rrto_{\alphabar_{f',g',h'}}& &
(f' \tenAbar g') \tenAbar h'.
\end{diagram} $$
The simple observation that $\phi \tenAbar( \gamma \tenAbar \eta ) =
V_{i,j,k}(\phi,\gamma,\eta,Id_0)$ with a similar identity for the
right vertical map and the fact that $V_{i,j,k}$ is a functor shows
that the square indeed commutes.

We emphasize that the reassociation $\alphabar$ is not
coherent.

\vspace{1 cm}

There is a well-known categorical principle that in any
category
$ \bf C$, $Hom_{\bf C}(X,G)$ inherits the structure
of the codomain object. If a category $\bf B$ has some higher-order
structure (we remain suitably ambiguous) then, loosely speaking,
$$ \coprod\limits_{n \in {\bf N^+}} {\bf Hom_{Cat}(X}^n {\bf ,B)}$$
inherits this structure. In our case we set $\bf X= \bf B$.
We claim that if $\bf B$ has the structure of an associative category
then so does
$$\Bbar = \coprod\limits_{n \in {\bf N^+}} \homBnB . $$
Since this construction will be used often, we feel obliged to go
through the gory details at least once. The tensor and reassociation
are defined as follows.
Let $f:{\bf B^q} \longrightarrow \bf B$, $g:{\bf B^r}
\longrightarrow \bf B$ and $h:{\bf B^s} \longrightarrow \bf B$. Then $$
\begin{diagram}
f \tenBbar g = \otimes_{\bf B} \circ (f \times g) : {\bf B^q
\times B^r}
\rrto^<{f \times g} & &{\bf B \times B} \rrto^{\otimes_{\bf B}} & &
\bf B.\end{diagram} $$
The reassociation
$$
\begin{diagram}
      & & \\
{\bf B^{q+r+s}} \rrtod_{(f \tenBbar g) \tenBbar h} \rrtou^{f \tenBbar (
g
\tenBbar h)}  &\Downarrow \betabar_{f,g,h} & {\bf B}\\
 &  &
\end{diagram}$$
is set to $$ \betabar_{f,g,h}(b_1, b_2 , \ldots ,b_{q+r+s})=
\beta_{f(b_1, \ldots, b_q),g(b_{q+1}, \ldots, b_{q+r}),h(b_{q+r+1},
\ldots
, b_{q+r+s})}. $$

We have proven the following.
\begin{lem}
Given an associative category $(\bf B ,\otimes_{\bf B}, \beta_{\bf B})$,
then $(\Bbar, \tenBbar,\betabar)$ also has the an associative category
structure.
\end{lem}
\begin{defi}
Given an associative category $(\bf B ,\otimes_{\bf B}, \beta_{\bf B})$,
the category of iterates of $\otimes_{\bf B}$, denoted  ${\overline
{\bf It(B,
\otimes_B,
\beta_B)}}$ or $\Bit$, is the associative subcategory of
$(\Bbar,
\tenBbar,\betabar)$ generated by $Id_B \in {\bf Hom(B,B)}$.
\end{defi}
$\Bit$ can be looked at as a disjoint union
$ \coprod\limits_{n \in {\bf N^+}} \Bitn$ where $\Bitn$ is constructed
recursively. ${\bf It(B)_1}=Id_{\bf B}.$ $\Bitn$ has as objects
$\otimes_{\bf B} \circ (f \times g)$ where $f$ is an object in
${\bf It(
B)_a}$, $g$ is an object in ${\bf It(B)_b}$ and $a+b=n$. Morphisms
of
$\Bitn$ are generated like the objects, however, there are also
morphisms inherited from the associative category structure of $\Bbar$.

\begin{defi}[Strict Tensor Functor]
Given two associative categories $(\bf B, \otimes,\beta)$ and $({\bf
B'},
\otimes ',\beta ')$, a strict tensor functor between them is a
functor $F: \bf B \longrightarrow {\bf B'}$ satisfying the following
two requirements:

i) $F(A \otimes B) = F(A) \otimes ' F(B) $ and

ii)$F(\alpha_{A,B,C}) = \alpha'_{FA,FB,FC}$
\end{defi}

\begin{prop}[Universality of $\Abar$]
For every associative category $(\bf B, \otimes,\beta)$, there is a
unique strict tensor functor
$$ P^{\bf B}: (\Abar = \coprod \An) \longrightarrow (\Bit = \coprod
\Bitn)$$ such that

i) $* \in {\bf A_1} \longmapsto Id_{\bf B} \in {\bf It(B)_1}$

ii) $* \in {\bf A_2} \longmapsto \otimes \in {\bf It(B)_2}$

iii) $\alpha$ (the isomorphism in ${\bf A_3}$) $\quad \longmapsto \quad
\beta$ (the isomorphism in ${\bf It(B)_3}$)
\end{prop}
{\bf Proof.} Any
functor $S: \Abar \longrightarrow \Bit$
has as its source a coproduct and hence can be described by its
components
$S_n:
\An
\longrightarrow \Bit$. Since we require  the first summand ${\bf
A_1}$
to go into the first summand, ${\bf It(B)_1}$, and S is a strict
tensor
functor, a quick induction shows that each $S_n$ actually lands in the
$n$-th summand of $ \Bit $ i.e. $\Bitn$. By the requirements of the
theorem
$P_1$,
$P_2$ and $P_3$ are forced. We construct and show uniqueness of $P_n$
with
the following argument. Assume $P_1, P_2, \ldots , P_{n-1}$ are
defined. Then $P_n$ is the unique functor making the
following pushout complete:
$$
\begin{diagram}
\coprod\limits_{i+j+k = n} {\bf A_i \times A_j \times A_k \times
\dot{I}}
\ddto|<<\ahook
\rrto^{W_{i,j,k}} & &
\coprod\limits_{a+b = n} {\bf A_a \times A_b} \ddto^{U_{a,b}}
\xto[3,1]^{\varphi_{a,b}}\\
\\
\coprod\limits_{i+j+k = n} {\bf A_i \times A_j \times A_k \times I}
\rrto^{V_{i,j,k}} \xto[1,3]^{\psi_{i,j,k}} & &{\bf
A_n}\ar@{-->}[dr]^{P_n}\\
 & & &\Bitn .
\end{diagram}
$$
Where $$\varphi_{a,b}(f,g)=P_a(f) \tenBbar P_b(g) =
P_2(*)(P_a(f) \times P_b(g)$$
$$\psi_{i,j,k}(f,g,h,0) = P_i(f) \tenBbar (P_j(g) \tenBbar P_k(h) )=
P_3(0)(P_i(f) \times P_j(g) \times P_k(h))$$
$$\psi_{i,j,k}(f,g,h,1) = (P_i(f) \tenBbar P_j(g)) \tenBbar P_k(h) =
P_3(1)(P_i(f) \times P_j(g) \times P_k(h)).$$
Let $\phi : f \longrightarrow f'$ , $\gamma: g \longrightarrow g'$ and
$\eta : h \longrightarrow h'$ then
\begin{eqnarray*}
\psi_{i,j,k}(\phi, \gamma, \eta ,\iota)  &=&P_3(Id_1)(P_i(f) \times P_j(g) \times
P_k(h)) \circ \betabar_{f,g,h}\\
    &=& \betabar_{f',g',h'} \circ P_3(Id_0)(P_i(f) \times P_j(g) \times
P_k(h))
\end{eqnarray*}

A simple diagram trace shows that the outer square commutes.
The pushout
insures that there is a unique map $P_n : \An \longrightarrow \Bitn $.

In order to show that $P$ is a strict tensor, let $f$ and $g$ be objects
in ${\bf A_a}$ and ${\bf A_b}$ respectively. The following square
commutes.
$$ \begin{diagram}
f,g \ar@{|->}[rrrr] \ar@{|->}[dddd] & & & & U_{a,b}(f,g)
\ar@{|->}[dddd] \\
            &\Abar \times \Abar \rrto^\tenAbar \ddto_{P \times P}& &
\Abar \ddto^P
\\
\\
            &\Bbar \times \Bbar \rrto_\tenBbar & & \Bit \\
P_a(f),P_b(g) \ar@{|->}[rrrr] & & & &\Text{$P_n(U_{a,b}(f,g))=$ \\
$P_a(f)
\tenBbar P_b(g)$}.
\end{diagram} $$
The lower right-hand equality holds because that is just the upper
right-hand triangle in the previous pushout diagram. Similar
arguments are needed for morphisms. $P$'s uniqueness follows from
the fact that $\varphi$ and $\psi$ are used in the proof of $P$
being strict.$\Box$

\section{$\bf A$, The 2-sketch of associative categories}

In order to define the 2-sketch, we need to give the $\An$ an operad
structure. An operad is a way of describing the structure of all the operations
that an algebraic object has. There are many different definitions of 
operads, but we shall use the simplest definition that will meet 
our needs.  We
basically follow
May's \cite{May} definition of a topological $A_\infty$ operad.
For a general introduction to operads, see the first few papers
of \cite{Loday}. For a very general definition of an operad, see
\cite{KrizMay}.
We will give here the definition of a non-$\Sigma$ non-unital operad in ${\bf Cat}$.
\begin{defi}An {\em operad} $\OP$ in ${\bf Cat}$ 
consists of categories  $\OP_j, j\ge 1$
and composition functors
$$ Q:\OP_k \times \OP_{j_1} \times \cdots \times \OP_{j_k}
\longrightarrow \OP_{j_1+j_2+ \cdots +j_k}$$
satisfying the following condition. If $\Sigma j_s = j$, $\Sigma i_t=i$,
 $g_s=j_1+j_2+ \cdots +j_s$, and $h_s=i_{g_{s-1}+1} + \cdots + i_{g_s}$
for $1 \le s \le k$, then the following ``associativity'' diagram 
commutes:
\end{defi}
\begin{equation}
 $$\xymatrix{   
 \OP_k \times \prod_{s=1}^k \OP_{j_s} \times \prod_{r=1}^j \OP_{i_r}
 \ar[rrrr]^{Q \times Id}
 \ar[dddd]_{Shuffle}
 & & & &
 \OP_j \times \prod_{r=1}^j \OP_{i_r}
 \ar[dd]^Q
\\
\\
 & & & &
\OP_i
\\
\\ 
 \OP_k \times (\prod_{s=1}^k \OP_{j_s} \times (\prod_{q=1}^{j_s}
 \OP_{i_{g_{s-1}+q})))}
\ar[rrrr]_{Id \times \prod_s Q}
&&&&
\OP_k \times \prod_{s=1}^k \OP_{h_s}.
\ar[uu]_Q
  }$$
\end{equation}

Each $\OP_j$ is to be thought of as the structure of the $j$-ary operations.
 $Q(c,d_1, \ldots,d_k)$ is the composition of the $c$ operation with the 
product of the $d_s$.

An operad in our paper is 
a way of combining associations. Given an association of $n$ letters and
$n$ associations of $m_1, m_2, \ldots , m_n$ letters, an operad makes
a new association of $t=\sum m_i$ letters by considering each of the
$m_i$ letters to be one unit of the original association.

In order to  make the discussion of operad  more readable, we conform to the 
following conventions.
When we refer to partitions of a set X, we mean an
{\it
ordered} set of disjoint nonempty subsets of X.
Every (ordered) partition of $t$ objects into $n \le t$ disjoint subsets
can be written
as an order-preserving surjection $\pi:t \longrightarrow n$.

\begin{prop}
${\bf A_1, A_2, \ldots , A_n, \ldots}$ has the structure of an operad.
\end{prop}
{\bf Proof.} We define
  a functor
$$Q_{n,t,\pi}:{\bf A_n \times A_{m_1} \times A_{m_2} \times \cdots
\times
A_{m_n}} \longrightarrow {\bf A_t} $$
for all $n$, for all $t \ge n$  and for all partitions $\pi: t
\longrightarrow
n$ where $m_i = | \pi ^{-1}(i) |$. We define the $Q_{n,t,\pi}$'s by
induction on n. For $n=1$, and for any $t \ge 1$ there is a unique $\pi
: t
\longrightarrow 1$. We set $$ Q_{1,t,\pi} = Proj: {\bf A_1 \times A_t}
\longrightarrow {\bf A_t}.$$
Now assume every $Q_1, Q_2, \ldots , Q_{n-1}$ is defined for each
$t$ and partition $\pi$. We define $Q_{n,t,\pi}$ for a
given $t$ and $\pi$ on 0)objects and 1)morphisms.

0) Let $f$ be an object in $\An$ and $g_1, g_2, \ldots, g_n$ be objects
in ${\bf A_{m_1}, \ldots , A_{m_n}}$ respectively. Since each object $f$
in $\An$ has the property that $f=U_{a,b}(f_1,f_2)   ~~(a+b = n)$ for
unique $a,b,f_1$ and $f_2$. We let
\begin{equation} Q_{n,t,\pi}(f,g_1,g_2, \ldots , g_n ) =
\label{defofQ1} \end{equation}
$$U_{a,b}(Q_{a,t_1,\pi_1}(f_1,g_1, g_2, \ldots , g_a), Q_{b,t_2, \pi_2}
(f_2, g_{a+1}, g_{a+2}, \ldots , g_{a+b})).$$
where
$$ t_1 = \sum\limits_{i=1}^a m_i = \sum\limits_{i=1}^a |\pi^{-1}(i)| ,
\qquad
 t_2 = \sum\limits_{i=a+1}^{a+b} m_i = \sum\limits_{i=a+1}^{a+b}
|\pi^{-1}(i)|
$$ and $\pi_1$ amd $\pi_2$ are the restrictions of $\pi$ to $a$ and $b$.

1)Since $Q_{n,t,\pi}$ is a functor, it suffices to define it on the
generating morphisms of $\An$.

1a)Morphisms in $\An$ of the form $\phi=U_{a,b}(\phi_1,\phi_2)$ are done
similarly as (2).

1b)Morphisms in $\An$ of the form $\phi =
V_{i,j,k}(\phi_1,\phi_2,\phi_3,\iota)$ are handled in the following
manner:
Let
$$ t_i = \sum\limits_{x=1}^i m_x , \qquad
t_j = \sum\limits_{x=i+1}^{i+j} m_x , \qquad
t_k = \sum\limits_{x=i+j+1}^{i+j+k} m_x . $$
Then \begin{equation}Q_{n,t,\pi}(\phi,\gamma_1,\gamma_2, \ldots,
\gamma_n)
= \label{defofQ2} \end{equation} 
$$
V_{t_i,t_j,t_k}(Q_{i,t_i,\pi_i}(\phi_1,\gamma_1,\ldots,\gamma_i),
Q_{j,t_j,\pi_j}(\phi_2,\gamma_{i+1},\ldots,\gamma_{i+j}),
Q_{k,t_k,\pi_k}(\phi_3,\gamma_{i+j+1},\ldots,\gamma_{n}),\iota )
$$

An example of the way the operad works is
in order. Let
$n=4$ and
$t=13$. Let $f \in {\bf A_4}$ correspond to the following association:
$A[[BC]D]$.
Let $g_1, g_2, g_3$ and $g_4$
correspond to the following associations:$$A(BC), \quad D, \quad
(EF)((GH)I), \quad ((JK)L)M. $$ Then $Q(f,g_1,g_2,g_3,g_4)$
corresponds to
$$A(BC)[[ D  ((EF)((GH)I))] ((JK)L)M]. $$

For an example of the way reassociations are handled, let
$\phi$
correspond to the reassociation $A[[BC]D] \longrightarrow [AB][CD]$ and
$$\begin{diagram}
A(BC) \ddto^{\gamma_1} & D\ddto^{\gamma_2}  & (EF)((GH)I)\ddto^{\gamma_3}  & (JK)L)M \ddto^{\gamma_4}
\\
\\
(AB)C & D & ((EF)(GH))I & J(K(LM)).
\end{diagram} $$
Then the operad will produce: $$\begin{diagram}
A(BC)[[ D  ((EF)((GH)I))]  (JK)L)M ] \ddto^{Q(\phi,\gamma_1,\gamma_2,
\gamma_3, \gamma_4)}
\\
\\
[((AB)C) D ] [ (((EF)(GH))I) ( J(K(LM)))].
\end{diagram}$$

If the $\An$ is to be a real operad, the $Q_{n,t,\pi}$s must satisfy the
``associativity'' diagram of Definition 4. We abandon all
unnecessary
subscripts for the benefit of the reader. $\vec{f}$ is used to mean a
sequence of
$f$'s of the ``right'' length. The lemma for morphisms is left to the
reader.
\begin{lem}[Associativity of Q]
Let $h$ be an object in $\An$. Let $g_i$ be objects ${\bf A_{m_i}}$.
We have the following equality:
$$ Q(h,Q(g_1,\vec{f}),\ldots,Q(g_n,\vec{f})) =
Q(Q(h,g_1,\ldots,g_n),\vec{f})$$
\label{assocofQ} \end{lem}
{\bf Proof.} By induction on $n$. If $n=1$ then $h=* \in {\bf A_1}$ and
$$LHS=Q(*,Q(g_1,f,\ldots,f))=
Q(g_1,f,\ldots,f)=Q(Q(*,g_1),f,\ldots,f)=RHS.$$
Assume the lemma is true for all $a < n$. Then $h\in \An$ and
$h=U_{a,b}(h_a,h_b)$ for unique $a,b,h_a$ and $h_b$. We have
\begin{eqnarray*}
LHS &=_1& Q(h,Q(g_1,\vec{f}),\ldots,Q(g_n,\vec{f})) \\
    &=_2& Q(U_{a,b}(h_a,h_b),Q(g_1,\vec{f}),\ldots,Q(g_n,\vec{f})) \\
    &=_3&U_{a,b}(Q_a(h_a,Q(g_1,\vec{f}),\ldots,Q(g_a,\vec{f})),
               Q_b(h_b,Q(g_{a+1},\vec{f}),\ldots,Q(g_{a+b},\vec{f}))) \\
    &=_4&U_{a,b}(Q(Q_a(h_a,g_1,\ldots,g_a),\vec{f}),
               Q(Q_b(h_b,g_{a+1},\ldots,g_{a+b}),\vec{f})) \\
&=_5&Q(U(Q_a(h_a,g_1,\ldots,g_a),Q_b(h_b,g_{a+1},\ldots,g_{a+b})),\vec{f})\\
             &=_6& Q(Q(U_{a,b}(h_a,h_b),g_1,\ldots,g_{a+b}),\vec{f}) \\
             &=_7& Q(Q(h,g_1,\ldots,g_n),\vec{f}) \\
             &=& RHS.
\end{eqnarray*}
$=_2$ and $=_7$ are from the definition of $h$. $=_3$, $=_5$ and $=_6$
are from definition of $Q$. $=_4$ is from the induction
hypothesis.$\Box$

A 2-sketch (called an algebraic 2-sketch in \cite{Heller}) is a strict
tensor
2-category whose underlying category
(0-cells and 1-cells) is a sketch i.e. a sketch that is enriched over
{\bf Cat}. An algebra $F$ for a 2-sketch $\bf G$ is a strict tensor
2-functor $F:{\bf G} \longrightarrow {\bf Cat}$.

At last, we are ready to define the 2-sketch, $\bf A$, of the theory of
associative categories. $\bf A$ is a strict tensor 2-category. The
objects are the positive natural numbers. In order for $\bf A$ to be a
2-category, it must be enriched over {\bf Cat} i.e. every hom set must
be a category and composition must be a functor. Given any two positive
integers $n$ and $k$, we define the category
$$\homnk =\coprod\limits_{\pi :n \rightarrow k} {\bf A_{|\pi^{-1}(1)|}
\times A_{|\pi^{-1}(2)|} \times \cdots \times A_{|\pi^{-1}(k)|} } =
 \coprod\limits_{\pi :n \rightarrow k}
\prod\limits_{i=1}^k {\bf A_{|\pi^{-1}(i)|}}$$
where $\pi$ ranges over all partitions of $n$ into $k$ parts. Notice
that $$\homnk = \left\{
\begin{array}{r@{\quad:\quad}l}
\An & k=1 \\
{\bf (A_1)^k} = {\bf A_1} & k=n \\
{\bf \emptyset}& k > n
\end{array} \right. .$$
Each object of $\homnk$ corresponds to a partial association of $n$
letters into $k$ parts. Each of the $k$ parts is totally associated. For
example, a typical object in ${\bf Hom(10,4)}$ looks like this
$$(AB), \quad (C(DE))F, \quad G, \quad H(IJ).$$ This object corresponds
to the partition
$10 = 2 + 4 + 1 + 3$. We write objects of $\homnk$ as
$$f=(f_1,f_2, \ldots ,f_k)$$ where each $f_i$ is an object of ${\bf
A_{|\pi^{-1}(i)|}}$.
Morphisms of $\homnk$ correspond to reassociations of partial
associations. Since $\homnk$ is made up of a disjoint union, there are
only reassociations of partial associations {\it of the same partition}.
Morphisms are written as
$$ \phi = (\phi_1, \phi_2, \ldots, \phi_k). $$

Composition:
Consider the following situation:
$$\begin{diagram}
n \rrto^{f'}\rrtou^f \rrtod_{f''}&\Text{$\Downarrow \phi$\\ \\
$\Downarrow
\phi
'$}
& k
\rrto^{g'}\rrtou^g \rrtod_{g''}&\Text{$\Downarrow \gamma$\\ \\
$\Downarrow
\gamma
'$}
& l.
\end{diagram}$$
Let $f=(f_1,f_2, \ldots , f_k)$ with its corresponding partition
$\pi_f:n \longrightarrow k$. Let $g=(g_1,g_2,\ldots, g_l)$ with its
partition $\pi_g:k \longrightarrow \l$. Then horizontal composition,
$\circ_H$ , is defined as
$$ g \circ_H f = h = (h_1,h_2,\ldots,h_l)$$
with corresponding partition $\pi_h = \pi_g \circ \pi_f:n
\longrightarrow l $, where $$ h_i = Q_{|\pi^{-1}_g(i)|,|\pi^{-1}_h(i)|,
\pi_f|}(g_i,f_{i_1},f_{i_2}, \ldots , f_{i_s}). $$
The $i_j$ range over $\pi^{-1}_h(i)$ and $\pi_f|$ is a restriction of
$\pi_f$ to this subset. An example is called for. Let $n=10,~~~ k=4$ and
$l=2$. Let $f$ correspond to $$(A(BC)), \quad D, \quad (E(F(GH))), \quad
(IJ). $$ Let $g$ correspond to $$[AB], \quad [CD]. $$ Then $g \circ_H f$
corresponds to $$[(A(BC)) D], \qquad [(E(F(GH))) (IJ)]. $$

Horizontal composition of 2-cells is also done with the operad $Q$.
Associativity of horizontal
composition follows
from lemma \ref{assocofQ}. We leave the details to the reader, however,
it is obvious once you look at what the $Q$'s are constructed to do.

Vertical composition of 2-cells, $\circ_V$, is much simpler. Let
$\phi=(\phi_1,\phi_2, \ldots, \phi_k)$ be of a reassociation of a
particular partition
and let $\phi'=(\phi_1',\phi_2', \ldots, \phi_k')$ be of the {\it same }
partition, then
$$\phi' \circ_V \phi = (\phi'_1 \circ \phi_1, \phi'_2 \circ \phi_2,
\ldots, \phi'_k \circ \phi_k).$$ Associativity of vertical composition
is obvious.

In order for $\bf A$ to be an honest 2-category, horizontal and vertical
composition must ``commute'' i.e.
\begin{lem}
$$(\gamma' \circ_H \phi') \circ_V (\gamma \circ_H \phi) = (\gamma'
\circ_V \gamma) \circ_H (\phi' \circ_V \phi) :g \circ_H f
\longrightarrow g' \circ_H f'. $$
\end{lem}
{\bf Proof.} LHS=$(\lambda_1, \lambda_2, \ldots, \lambda_l)$ where
$$\lambda_i = Q_{|\pi^{-1}_{\gamma'}(i)|, m, \pi_{\phi'}|}(\gamma'_i,
\phi'_{i_1}, \phi'_{i_2}, \ldots , \phi'_{i_m}) \circ
Q_{|\pi^{-1}_{\gamma}(i)|, m,
\pi_{\phi}|}(\gamma_i,\phi_{i_1}, \phi_{i_2}, \ldots , \phi_{i_m}).$$
RHS=$(\rho_1, \rho_2, \ldots, \rho_l)$ where
$$\rho_i = Q_{|\pi^{-1}_{\gamma}(i)|, m, \pi_{\phi}|}(\gamma'_i \circ
\gamma_i,\phi'_{i_1} \circ \phi_{i_1}, \phi'_{i_2} \circ
\phi_{i_2},\ldots,\phi'_{i_m} \circ \phi_{i_m}).$$
Since $\pi_{\gamma} = \pi_{\gamma'}$ and $\pi_{\phi} = \pi_{\phi'}$,
these
two $Q$'s are actually the same functor and by the functoriality of $Q$,
 $\lambda_i = \rho_i$ and hence LHS = RHS. $\Box$
\\
\begin{prop}
$\bf A$ has a strict 2-tensor structure, $\tenA$.
\end{prop}
{\bf Proof.}

0) 0-cells: $n \tenA m = n + m$.

1) 1-cells:
Let $f=(f_1,f_2, \ldots , f_n)$ and $g=(g_1,g_2,\ldots, g_m)$ be two
1-cells, then $f \tenA g = (f_1,f_2, \ldots , f_n, g_1,g_2,\ldots,
g_m)$.

2) 2-cells: This is done the same way as 1-cells. $\Box$

In order for the following definition to make sense, we must remind
ourselves that ${\bf Cat}$, the 2-category of categories, functors and
natural transformations, has a strict 2-tensor structure: product. (We
parenthetically note that product is only coherently associative but we
think of it - perhaps in error - as strict because of the usual
coherence theories.)

\begin{defi}
Let $\homAcat$ be the category (we forget, for now, its higher-order
structure) whose objects are the strict tensor 2-functors $$R: \bf A
\longrightarrow {\bf Cat}$$
i.e. the 2-functors that satisfy

0) $R(n \tenA m) = R(n) \times R(m)$,

1) $R(f \tenA g) = R(f) \times R(g)$ and

2) $R(\phi \tenA \gamma) = R(\phi) \times R(\gamma)$.
\\
\\
Morphisms are strict 2-natural transformations $$F:R \Longrightarrow
S$$ i.e.
0) For every 0-cell $n$ in $\bf A$  there is a functor $F(n):R(n)
\longrightarrow S(n)$ such that $F(n \tenA m) = F(n) \times F(m)$.

1) For every 1-cell $f:n \longrightarrow m$ in $\bf A$ the following
square commutes ``on the nose''
$$\begin{diagram}
R(n) \rrto^{F(n)} \ddto_{R(f)} & & S(n) \ddto^{S(f)}\\
\\
R(m) \rrto_{F(m)} & & S(m). \end{diagram} $$

2) For every 2-cell $\phi: f \Longrightarrow f'$ in $\bf A$, the
following square commutes ``on the nose''
$$\xymatrix{
R(n) \ar[rrrr]^{F(n)} \ar@/^/[ddddd]^{R(f)} \ar@/_/[ddddd]_{R(f')} & & & & S(n)
\ar@/^/[ddddd]^{S(f)}
\ar@/_/[ddddd]_{S(f')}
\\ \\
\Rightarrow & &  & & \Rightarrow 
\\
\\ \\
R(m) \ar[rrrr]_{F(m)} & && & S(m) 
}$$ where the left and right 2-cells are $R(\phi)$ and $S(\phi)$.
\end{defi}

\begin{defi}
The category ${\bf Assoc-Cat_{st}}$ has as objects associative
categories and as morphisms, strict tensor functors.
\end{defi}

\begin{prop}[$\bf A$ as 2-sketch of associative categories]
The category $\homAcat$ is equivalent to ${\bf Assoc-Cat_{st}}$.
\end{prop}
{\bf Sketch of Proof.} The proof calls for only a few minutes of
staring at the definitions.
We will not go through all the hideous details; but shall point
the way. Given a strict tensor 2-functor $R:\bf A \longrightarrow {\bf
Cat}$ we set the underlying category, $\bf B$, to be $R(1)$. Set
$\otimes_{\bf B} = R(f)$ where $f$ is the unique morphism from $2$ to
$1$ in $\bf A$. $\beta_{\bf B} = R(\iota)$ where $\iota$ is the unique
nontrivial isomorphism in ${\bf A_3}$.

To a strict 2-natural transformation $F:R \Longrightarrow S$ we assign a
strict tensor functor (also called $F$) $F:R(1) \longrightarrow S(1)$.
To every positive natural number $n=1+1+ \cdots +1$, there is a functor
$$[F(n)=F(1)^n] : [R(n)=
R(1)^n] \longrightarrow
[S(n) = S(1)^n].$$
If we let $R(1)= \bf B$ and $S(1) = {\bf B'}$, then the above line looks
like the more familiar $$F^n:{\bf B^n} \longrightarrow {\bf B'^n}.$$ The
two commuting diagrams in the definition of strict 2-natural
transformations correspond to the two requirements for a functor to be
a strict tensor functor. $\Box$

\section{Associative categories with units}
\begin{defi}
An associative category with a unit is an associative category $(\bf B,
\otimes_{\bf B}, \beta)$ with a distinguished object $I \in \bf B$ and
the following two natural isomorphisms:
$$ L_A:I \otimes A \longrightarrow A \qquad R_A: A \otimes I
\longrightarrow A.$$
\end{defi}
There are times when the following coherence condition will be
important.
An associative category with a unit in which the following diagram
commutes is said to have {\it unital coherence}:
$$\xymatrix{
(A \otimes I) \otimes B \ar[rrrr]^{\beta_{A,I,B}} \ar[rrdd]_{R_A \otimes
Id_B} &&&& A \otimes (I \otimes B) \ar[lldd]^{Id_A \otimes L_B} \\
\\
&& A\otimes B }$$
If the $L_A$ and $R_A$ are identity natural isomorphisms then we say
the associative category is {\it unital strict}.

The $L_A$'s and $R_A$'s are isomorphisms connecting associations of
$n+1$ letters to associations of $n$ letters. In the free associative
category, they would be isomorphisms from the objects of ${\bf
A_{n+1}}$ to the objects of $\An$. These
isomorphisms can be formalized with a new inductive scheme of
pushouts. We will only sketch this here.
In this new formalism, we generate a
new sequence of groupoids ${\bf A'_0 , A'_1, \ldots , A'_n , \ldots}$.
Notice that this time we have a ${\bf A'_0}$ whereas there is no
${\bf A_0}$. Each ${\bf A'_n}$ is actually a subgroupoid of
${\bf A'_{n+1}}$. The free associative category with unit will not be
the coproduct of the ${\bf A'_n}$ rather it will be the colimit i.e.
$${\bf {\bar A'}} = Colim_{n \ge 0} {\bf A'_n}. $$
The details of the structure of ${\bf {\bar A'}}$ and its universal
properties are straightforward.

The ${\bf A'_n}$ also have the structure of an operad and this operad is
used to construct ${\bf A'}$, the 2-sketch of associative categories
with units. There is one interesting difference that is worth pointing
out. For ${\bf A}$ we had
$$\homnk =
\coprod\limits_{\pi :n \rightarrow k}
\prod\limits_{i=1}^k {\bf A_{|\pi^{-1}(i)|}}.$$
where the
partitions
$\pi$'s are surjective. Here --- when talking about units --- we allow
nonsurjective $\pi$'s. $|\pi^{-1}(i)|$ can
equal $0$ and we would get ${\bf A'_0}$ which corresponds to the unit.

\section{Maximal indiscrete subgroupoids}

The rest of this paper is dedicated to calculating the
fundamental group of each of the
$\An$ and their quotients. If the fundamental group of a quotient of
$\An$ is trivial, then there is only one path between any two
vertices
(objects) in the complex (groupoid) and it is called coherent.

In order to determine the fundamental group of a simplicial complex, one
can use the method of maximal trees
(see e.g. \cite{Rotman}
or \cite{Lynd&Sch}.)
Given any complex $K$, in order to find the fundamental group of $K$,
denoted
$\pi(K)$,
one associates a maximal tree, $T_K$, to $K$. A tree in $K$ is a
connected subcomplex of $K$ which has no circuits. A maximal
tree in $K$ is a tree in $K$ contained in no larger tree or,
equivalently, a tree that contains all vertices of $K$. $\pi(K)$ is
then given by the following presentation
\begin{itemize}
\item Generators
\begin{enumerate}
\item All edges $(u,v)$ in $K$.
\end{enumerate}
\item Relations
\begin{enumerate}
\item $(u,v)=e$ if $(u,v)$ is in $T_K$.
\item $(u,v)(v,w)=(u,w)$ if $u,v,w$ lie in the same simplex of $K$.
\end{enumerate}
\end{itemize}
There is a standard theorem that $\pi(K)$ is (up to conjugation)
invariant under a change of the maximal tree.

When dealing with the Catalan groupoids, we employ a method analogous
to maximal
trees.  To
each groupoid $\An$, we will assign
$\Tn$ a maximal
indiscrete subgroupoid (henceforth MIS) - the categorical analog of a
maximal tree. A MIS is a subgroupoid
with the same objects and exactly one isomorphism between each
ordered pair of objects.
It must be stressed that $\Tn$ is an MIS and not a tree.
There may, in
fact, be circuits in our
MIS but they correspond to a commuting part of the groupoid.
(Both the language of trees and the
language of categories will be employed. ``vertex'' and ``object'' will
be interchanged, as will ``edge'' and ``morphism''.)

The MIS's are defined inductively. ${\bf T_1} = {\bf A_1} = \bf
1$. Assume ${\bf T_1,T_2, \ldots, T_{n-1}}$ are defined, then $\Tn$ is
constructed in a manner similar to pushout (1). We set $j=1$ in pushout
(1) and get the following pushout.
$$
\begin{diagram}
\coprod\limits_{i+1+k = n} {\bf T_i \times T_1 \times T_k \times
\dot{I}}
\ddto|<<\ahook
\rrto^{W_{i,1,k}|} & &
\coprod\limits_{a+b = n} {\bf T_a \times T_b} \ddto^{U_{a,b}|} \\
\\
\coprod\limits_{i+1+k = n} {\bf T_i \times T_1 \times T_k \times I}
\rrto^{V_{i,1,k}|} & &{\bf T_n}
\end{diagram}
$$
where $W_{i,1,k}|$ is the restriction of $W_{i,j,k}$ of our original
pushout. (This is intuitive but actually too swift because we have not
yet shown that $\Tn$ is a subcategory of $\An$. So define $W|$ in a
similar manner to the way $W$ was in pushout (1).)
$\Tn$ contains only one class of morphisms between each component. This
corresponds to moving the outermost parentheses one place to the
left. We use the fact that one can get from any association
of n letters to any other association of n letters by only moving the
parentheses one letter at a time.

Now inductively define $$L_n: \Tn \longrightarrow \An.$$ $L_1 =
Id_{\bf 1}:{\bf T_1} \longrightarrow {\bf A_1}.$ Assuming $L_1, L_2,
\ldots
, L_{n-1}$ are defined, $L_n$ is then constructed from the following
diagram:
\tiny
$$
\begin{diagram}
\coprod\limits_{i+j+k = n } {\bf A_i \times A_j \times A_k \times
\dot{I}}
\xto [3,0]|<<\ahook
\xto[0,3]^{W_{i,j,k}}& & &
\coprod\limits_{a+b = n} {\bf A_a \times A_b} \xto[3,0]^{U_{a,b}} \\
 & \coprod\limits_{i+1+k = n} {\bf T_i \times T_1 \times T_k \times
\dot{I}}
\dto|<<\ahook
\ulto^{\coprod \times L}
\rto &
\coprod\limits_{a+b = n} {\bf T_a \times T_b} \dto^{U_{a,b}|}
\urto^{\coprod L_a \times L_b}
\\
& \coprod\limits_{i+1+k = n} {\bf T_i \times T_1 \times T_k \times
I}
\dlto^{\coprod \times L}
\rto & {\bf T_n}\ar@{-->}[dr]^{L_n} \\
\coprod\limits_{i+j+k = n} {\bf A_i \times A_j \times A_k \times I}
\xto[0,3]^{V_{i,j,k}} &   & &{\bf A_n}
\end{diagram}
$$
\normalsize
The miters of the diagram are made of (co)products of $L_a$.  The upper
and
left-hand trapezoid commute because the miters are basically inclusions
and the parallel morphisms are defined the same way. A diagram chase
shows that $\An$ satisfies the inner pushout condition and hence there
is a unique map $L_n: \Tn \longrightarrow \An$.

In order for $\Tn$ to be a MIS of $\An$, $L_n$ must be
bijective on objects (maximal) and for every $t,t'$ in $\Tn$,
$Hom_{\Tn}(t,t')
=
*$,
the one object set (indiscrete).

$L_n$ can be shown to be bijective on objects with a short inductive
proof. The base case is true by definition. The inductive step follows
from the fact that the right hand trapezoid commutes;
$U_{a,b}$ is bijective on objects; and the product of
bijective-on-objects functors (the miters) is bijective-on-objects. So
going around
the right hand trapezoid are only functors that are
bijective-on-objects.

An inductive proof is used to show that $Hom(t,t') =*$.
Assume
${\bf T_a}$ and ${\bf T_b}$ are indiscrete subgroupoids.
Then the product of indiscrete subgroupoids are indiscrete
subgroupoids and hence ${\bf T_a
\times T_b}$ is an indiscrete subgroupoid. Since
there is only one class of morphisms
$(<i,1,k>)$ connecting these indiscrete subgroupoids, the
entire
$\Tn$ is indiscrete.

A few diagrams of the MIS are called for. We shall display the
generators of the groupoid of the first few $\Tn$ and the way they sit
in
$\An$. For each
$\An$ there are three types of generating morphisms:
\begin{enumerate}
\item those not in $\Tn$ - denoted
$$ \begin{diagram} \bullet \rrto & & \bullet \end{diagram}.$$
\item those in $\Tn$ within ${\bf T_a \times T_b}$ for some $a$ and $b$
- denoted
$$ \begin{diagram} \bullet \ar@2{->}[rr] & & \bullet \end{diagram}.$$
\item those in $\Tn$ of the form ${\bf A_i \times A_1 \times A_k \times
I} $ - denoted
$$ \begin{diagram} \bullet \ar@3{->}[rr] & & \bullet \end{diagram}.$$
\end{enumerate}
The first two maximal trees are simple ${\bf T_1} = {\bf A_1} = \bf 1$,
${\bf T_2} = {\bf A_2} = \bf 1$.
\begin{itemize}
\item  ${\bf T_3}$
$$ \begin{diagram} 
\Text{${\bf A_1 \times A_2} $\\$A(BC)$}
\ar@3{->}[rrr]^{\bf A_1 \times A_1 \times A_1 \times
I}_{\alpha}&  &  & \Text{${\bf A_2 \times A_1}	$\\$(AB)C$}
\end{diagram}$$
\item ${\bf T_4}$ $$ 
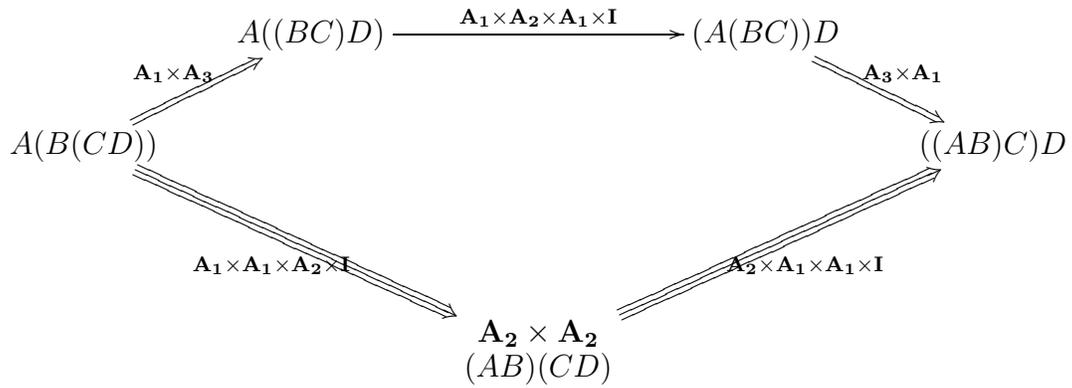
\begin{diagram}
			& A((BC)D) \rrto^{\bf A_1 \times A_2 \times A_1
\times I}&	 &(A(BC))D \ar@2{->}[dr]^{\bf A_3 \times A_1}\\
       A(B(CD))\ar@2{->}[ur]^{\bf A_1 \times A_3}
\ar@3{->}[ddrr]_{\bf A_1 \times A_1 \times A_2 \times I}& & &
&((AB)C)D\\
\\
  & &\Text{${\bf A_2 \times A_2}$\\$(AB)(CD)$} \ar@3{->}[uurr]_{\bf A_2
\times A_1 \times A_1
\times I}&
&
\end{diagram} $$
\item ${\bf T_5}$
{\setlength{\oddsidemargin}{0mm}
\tiny
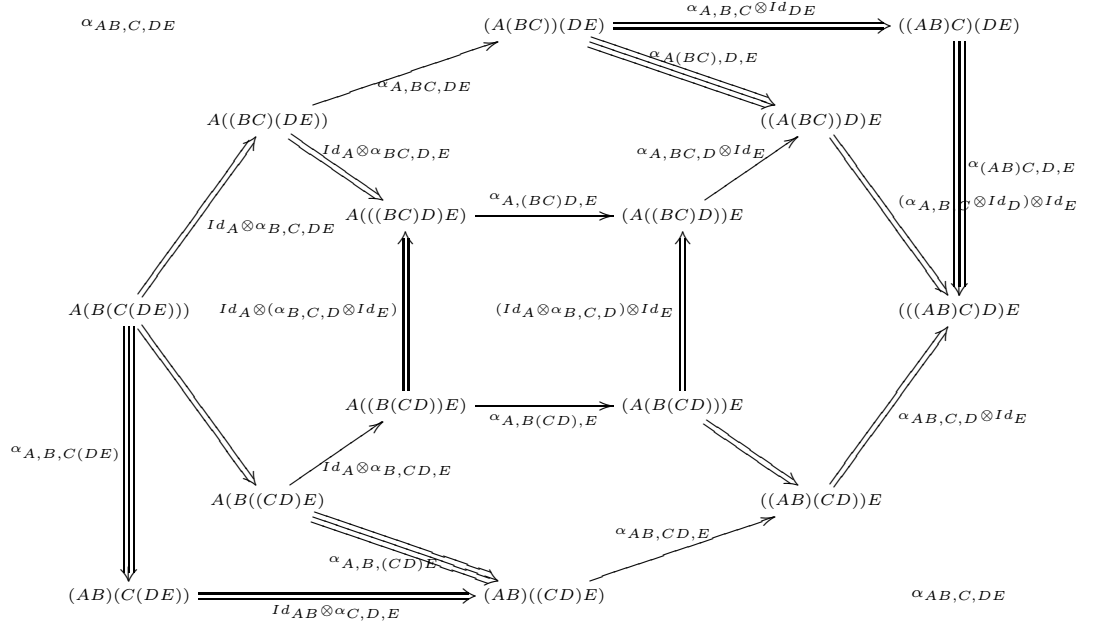
$$
\xymatrix@C=0pt{
\alpha_{AB,C,DE} &  &  &(A(BC))(DE)
\ar@2{->}[rrr]^{\alpha_{A,B,C} \otimes Id_{DE}}
\ar@3{->}[drr]^{\alpha_{A(BC),D,E}}
&  &  &((AB)C)(DE)
\ar@3{->}[ddd]^{\alpha_{(AB)C,D,E}}
\\
 &A((BC)(DE)) \urrto_{\alpha_{A,BC,DE}} \ar@2{->}[dr]^{Id_A \otimes
\alpha_{BC,D,E}}
&
&
&
&((A(BC))D)E
\ar@2{->}[ddr]^{(\alpha_{A,B,C} \otimes Id_D) \otimes Id_E}
&
\\
 &  &A(((BC)D)E) \rrto^{\alpha_{A,(BC)D,E}} &  &(A((BC)D))E
\urto^{\alpha_{A,BC,D} \otimes Id_E}
&
&
\\
A(B(C(DE))) \ar@2{->}[uur]_{Id_{A}\otimes \alpha_{B,C,DE}}
\ar@3{->}[ddd]_{\alpha_{A,B,C(DE)}}  \ar@2{->}[ddr] & & & & &
&(((AB)C)D)E
\\
 & &A((B(CD))E) \rrto_{\alpha_{A,B(CD),E}} \ar@2{->}[uu]^{Id_A
\otimes(\alpha_{B,C,D}\otimes Id_E)}
& &(A(B(CD)))E
\ar@2{->}[uu]^{(Id_A\otimes \alpha_{B,C,D})\otimes Id_E}  \ar@2{->}[dr]
\\
 &A(B((CD)E) \urto_{Id_{A}\otimes \alpha_{B,CD,E}}
\ar@3{->}[drr]_{\alpha_{A,B,(CD)E}}  &  &  &  &((AB)(CD))E
\ar@2{->}[uur]_{\alpha_{AB,C,D}\otimes Id_E} \\
 (AB)(C(DE)) \ar@2{->}[rrr]_{Id_{AB}\otimes \alpha_{C,D,E}}
& & &(AB)((CD)E)
\urrto^{\alpha_{AB,CD,E}}
& & & \alpha_{AB,C,DE}
}
$$}
\normalsize

Each vertex is reached by $\Tn$.
The only circuits are naturality squares. The long top map
is the same as the long bottom
map. The single-line arrows have more than one letter in the center.
The double-line arrows have only one letter in the center, but the map
is tensored i.e. old. The triple-line arrows have only one letter in the
center and are new i.e. not tensored.
\end{itemize}

{\bf Remark.} There is nothing canonical about our MIS. We could
have
chosen as our MIS morphisms those morphisms of the form
$<1,j,k>$, or some other scheme. This would have made a different
MIS but we would get --up to conjugation-- the same group at the end.

As with maximal trees, we must now collapse all the morphisms in
the MIS to a point. This is done with the following pushout in
the category of categories or groupoids:
$$ \begin{diagram}
\Tn \rrto^{!} \ddto_{L_n}
& & \bf 1 \ddto
\\ \\
\An \rrto
& &
\pi(\An,\Tn).
\end{diagram} $$
$\pi(\An,\Tn)$ can be thought of as the fundamental group of $\An$
relative to $\Tn$. Since we shall not change $\Tn$ in
this paper, we shorten $\pi(\An,\Tn)$ to $\pi(\An)$. The fact that
$L_n$ is surjective on objects, shows
that $\pi(\An)$ is a one object category. Since all the morphisms in
$\pi(\An)$ come from $\An$ or $\bf 1$ which only have
isomorphisms,
$\pi(\An)$ is a group. Every morphism in the MIS is sent to the
identity of $\pi(\An)$. The second type of relation comes
out of the pushout and the way $\Tn$ ``sits in'' $\An$.

\section{Presentation of the groups}
The generators of the groups are equivalence classes of generating
isomorphisms of the $\An$. Two generating isomorphisms are equivalent
if they are in the same image of $V_{i,j,k}$.
The image of
${\bf A_i \times A_j \times A_k \times I}$ under $V_{i,j,k}$, loosely
speaking, only contributes new morphisms. Each of these morphisms
are natural to one another i.e. they are two parallel sides of a
square that commutes under naturality. This can be seen by considering
the following situation.
Let $\phi :f \longrightarrow f'$ in ${\bf A_i}$,
      $\gamma:g \longrightarrow g'$ in ${\bf A_j}$
and $\eta : h \longrightarrow h'$ in ${\bf A_k}$. Then the following
square commutes out of the functoriality of $V_{i,j,k}$
$$ \begin{diagram}
U({\bf A_i \times A_{j+k}}) & & & U({\bf A_{i+j} \times A_k})\\
f \otimes (g \otimes h) \xto [0,3]^{V_{i,j,k}(Id_f,Id_g,Id_h,\iota)}
\ddto_
{V_{i,j,k}(\phi,\gamma,\eta,Id_0)}& & &
(f \otimes g) \otimes h \ddto^{V_{i,j,k}(\phi,\gamma,\eta,Id_1)}\\
\\
f' \otimes (g' \otimes h')
\xto [0,3]_{V_{i,j,k}(Id_{f'},Id_{g'},Id_{h'},\iota)}&
& &
(f' \otimes g') \otimes h'.
\end{diagram} $$
The left (respectively right) side of the diagram is in the image of
${\bf A_i \times A_{j+k}}$ under $U_{i,j+k}$ (resp. ${\bf A_{i+j}
\times A_k}$ under $U_{i+j,k}$.) The entire diagram is contained in
$\An$
where $n=i+j+k$. The point is that the entire image of $V_{i,j,k}$ is
really a set of edges of the $\An$ and of its $\Tn$. We shall
denote this set of edges as
$<i,j,k>$. So $<i,j,k>$ can be thought of as a set of morphisms $$
\begin{diagram} f \otimes (g \otimes h) \rrto
^{<i,j,k>} & &
(f \otimes g) \otimes h.
\end{diagram}  $$ one for each $f,g$ and $h$. We may denote the above
element of $<i,j,k>$ as $<i,j,k>_{f,g,h}$

The presentation of $\pi(\An)$ will be given by A set of
generators, $G_n$, and a set of relations, $R_n$.
We represent the operation of the
group as
$*$ and the trivial element as $e$.
The generators of the
group will be of the form $X<i,j,k>$ where X is an element of the free
monoid on two generators $\lam$ and $\rho$ (the monoid
operation is represented as
concatenation). Since the monoid is free, the length of an element is a
well defined concept. Intuitively new generators in $G_n$ are of the
form
$<i,j,k>$ where there is no prefix. Some
generators are from old components and the monoid is used to describe
which old
component the generator is from. If the generator is from
the right side of
${\bf A_a \times A_b}$ then we have $a$ blanks on the left
($\lambda^a$).
On the other hand, if the generator is on the left, we have $b$ blanks
to the right ($\rho^b$). In the latter case, for example, the generator
might have come even come from an earlier component and more $\lam$'s
and
$\rho$'s will be before the $<i,j,k>.$
For each generator $X<i,j,k>$ of $G_n$, the length of $X$ added to
the sum of $i,j$ and $k$ will equal $n$. In other words, a generator
first started in ${\bf A_{i+j+k}}$ and we use the free monoid to
describe how it sits in $\An$.
In order to facilitate writing out the generators, we define a
function $\#$. Given $X$, an element of the free monoid, and $G_a$, a
set of generators, we let
$$X\#G_a = \{ Xy:y \in G_a \}.$$ The relations are given as a set
of elements of the form $y_1=y_2$ by which we mean all elements
of the form
$y_1*y_2^{-1}$ are equal to $e$.
Given $X$ an element of the free monoid and $R_a$, a set of relations,
$$X\#R_a = \{ Xy_1*Xy_2* \cdots *Xy_m:y_1*y_2* \cdots *y_m \in R_a \}.$$

We begin by discussing the generators and relations and then give a
formal definition. The next section has a few calculations carried out.

\begin{itemize}
\item Generators.
The generators are equivalence classes of generating isomorphisms
(edges) in the $\An$. They form a set $G_n$. There are two types of such
morphisms:
\begin{enumerate}
\item Morphisms from the old components. This corresponds to the
upper right-hand corner of pushout (\ref{mainpo}). They are written as
$X<i,j,k>$ where $ln(X)+i+j+k = n$.
\item New morphisms between the old components. This corresponds to the
lower left-hand corner of pushout (\ref{mainpo}). They are written as
$<i,j,k>$ where $i+j+k = n$.
\end{enumerate}
\item Relations.
Relations come form:
\begin{enumerate}
\item commuting squares;
\item setting new morphisms that are in the MIS to the identity of
the group. i.e. $<i,j,k>=e$ if $j=1$;
\item old relations from the old components;
\item product of groupoids ${\bf
A_a
\times A_b}$. (as in topological spaces, generators of product
groupoids commute.)
\end{enumerate}
\end{itemize}

Some more words on the first type of
relation are needed.
The only commuting parts of $\An$ arise in the
following situation.
Let $\phi :f \longrightarrow f'$ in ${\bf A_i}$,
      $\gamma:g \longrightarrow g'$ in ${\bf A_j}$
$\eta : h \longrightarrow h'$ in ${\bf A_k}$ and
$$ \begin{diagram}
f \otimes (g \otimes h) \xto [0,2]^{<i,j,k>}
\dto_
{\phi\otimes(\gamma \otimes \eta)}& &
(f \otimes g) \otimes h \dto^{(\phi  \otimes
\gamma)\otimes\eta}\\
f' \otimes (g' \otimes h')
\xto [0,2]_{<i,j,k>}
& &
(f' \otimes g') \otimes h'.
\end{diagram} $$

This square commutes if the following diagram commutes:
$$ \begin{diagram}
f \otimes (g \otimes h) \xto[0,2]^{<i,j,k>}
\dto_{\phi\otimes(Id_g \otimes Id_h)} & &
(f \otimes g) \otimes h
\dto^{(\phi \otimes Id_g) \otimes Id_h}
\\
f' \otimes (g \otimes h)
\xto [0,2]^{<i,j,k>}
\ddto_{Id_{f'} \otimes (\gamma \otimes Id_h)} & &
(f' \otimes g) \otimes h
\ddto^{(Id_{f'} \otimes \gamma)\otimes Id_h}
\\ & (b) &
\\
f' \otimes (g' \otimes h)
\xto [0,2]^{<i,j,k>}
\dto_{Id_{f'}\otimes(Id_{g'} \otimes \eta)} & &
(f' \otimes g') \otimes h
\dto^{(Id_{f'} \otimes Id_{g'})\otimes\eta}
\\
f' \otimes (g' \otimes h')
\xto [0,2]^{<i,j,k>}  &  &
(f' \otimes g') \otimes h'.
\end{diagram} $$

Squares of these type generate all commuting
squares. Lets look at square $(b)$ in depth. The left vertical map,
$Id_{f'} \otimes (\gamma \otimes Id_h)$, is denoted by $\lam^i \rho^k y$
for some $y \in G_j$. Similarly, the right vertical map is denoted
$\rho^k \lam^i y$.
Since the square commutes
-- i.e. all four edges are ``in the same simplex''-- there is the
relation:
$$<i,j,k> * \rho^k \lam^i y = \lam^i \rho^k y * <i,j,k> .$$
The multiplication is written as regular multiplication rather then
morphism composition. Since $<i,j,k>$ denotes an {\em iso}morphism, in
fact, a whole set of isomorphisms, we can take inverses. Thus the above
equation looks like
$$<i,j,k> * \rho^k \lam^i y * <i,j,k>^{-1} = \lam^i \rho^k y.$$
This looks like the formula for HNN-extensions.
This is made more apparent by looking at the
categorical constructions of HNN-extensions. Given a group, $G$, and two
subgroups, $A,B$ with an isomorphism $f:A \longrightarrow B$ between
them, the HNN-extension, $H$, is given as the following pushout in the
category of groupoids:
$$ \begin{diagram}
A \times \dot{I} \rrto^W \ddto|<<\ahook & & G \ddto \\
\\
A \times I \rrto & & H
\end{diagram} $$
where $W(a,o)=a$ and $W(a,1)=f(a)$. Our pushout is then just a much more
complicated version of this. The $<i,j,k>$'s are to be thought of as the
new generators that extend all the old groups.

There is a special case of the above situation. When $j=1$ the $<i,j,k>$
generator is set to $e$ and the above equation looks like
$$ \rho^k \lam^i y = \lam^i \rho^k y.$$
This is the relation for a free product with amalgamation
of groups. This can be thought of as changing the $I$ in the lower
left-hand corner of the above pushout into $\bf 1$, the trivial one
object groupoid (with the left vertical map as the projection on $A$.)

We have only looked at the square $(b)$. There are, however, similar
equations for the other two types of boxes and they are given in the
scheme.

We give the inductive scheme for the generators and relations.
Throughout the scheme, $i,j$ and $k$ are positive integers.
$G_1 = R_1 = \emptyset$. Assume we have all the $G_k$ and $R_k$ for $k
\le n-1$ then $G_n$ and $R_n$ is given as:

\begin{itemize}
\item Generators. $G_n = $
\begin{enumerate}
\item $\{ [(\lam^i \# G_{n-i}) \cup
(\rho^{n-i} \# G_i) : i=1,2,\ldots, n-1\} \cup $
\item $\{<i,j,k> : i+j+k = n \}.$
\end{enumerate}
\item Relations. $R_n = $
\begin{enumerate}
\item For each $<i,j,k>$ such that $i+j+k = n$ we have the
unions of the following relations.
\begin{enumerate}
\item
$\{ <i,j,k> * \rho^k \rho^j x * <i,j,k>^{-1} = \rho^{j+k} x : x \in G_i
\}$
\item
$\{ <i,j,k> * \rho^k \lam^i y * <i,j,k>^{-1} = \lam^i \rho^k y : y
\in G_j \}$
\item
$\{ <i,j,k> * \lam^{i+j} z * <i,j,k>^{-1} = \lam^i \lam^j z : z \in
G_k \}$
\end{enumerate}
\item $\{<i,1,k>=e : i+1+k =n \} \cup$
\item $\{ [(\lam^i \# R_{n-i}) \cup
(\rho^{n-i} \# R_i) : i=1,2, \ldots , n-1 \} \cup $
\item $\{ \rho^{n-i} x * \lam^i y = \lam^i y * \rho^{n-i} x : x \in
G_i, y \in G_{n-i}
\} $
\end{enumerate}
\end{itemize}

\section{The first few groups}
In order to give the presentations in an clear fashion, we
conform to the following conventions about listing the generators and
relations: \begin{itemize}
\item The old edges are listed in columns below the names of the
components that contributed them. The new edges are listed at the
bottom.
\item If an old edge was set to $e$, we do not list it in further
groups. For example $<1,1,1>$ will be listed in $G_3$. Since it is set
to $e$, we do not list it or any of its ``progeny'' (e.g. $\lam <1,1,1>$
or $\rho \lam^3 \rho <1,1,1>$.)
\item We do not list relations about edges that were set to $e$.
\item Since the first nontrivial relation is in $R_6$ we do not list a
relation of type 3 until $R_7$.
\item Since the first nontrivial generator is in $G_4$, the first time
we have a relation of type 4 is in $R_8$. Due to the fact that $G_8$ and
$R_8$ will not be listed in this paper, we feel obliged to give this
relation of type 4:
$$\rho^4 <1,2,1> * \lam^4 <1,2,1> =\lam^4 <1,2,1> * \rho^4 <1,2,1>.$$
\item If two generators are set to be equal, then, in the future, we
only list one of them. We usually choose the shorter name e.g. $\lam^2
\rho <1,2,1>$ rather then $\lam \lam \rho <1,2,1>$.
\end{itemize}
Here are the groups:
$$\begin{array}{||c|c||}
\multicolumn{2}{c}{G_3}\\
\hline\hline
\bf A_1 \times A_2 & \bf A_2 \times A_1 \\
\hline
   & \\
\hline
\multicolumn{2}{||c||}{<1,1,1>}\\
\hline\hline
\end{array}$$

$$\begin{array}{||c||}
\multicolumn{1}{c}{R_3} \\
\hline\hline
<1,1,1> = e \\
\hline\hline
\end{array}$$
So $\pi({\bf A_3})$ is the trivial group.
\vspace{2cm}

$$\begin{array}{||c|c|c||}
\multicolumn{3}{c}{G_4}\\
\hline\hline
\bf A_1 \times A_3 & \bf A_2 \times A_2 & \bf A_3 \times A_1\\
\hline
\lam <1,1,1>& & \rho <1,1,1> \\
\hline
\multicolumn{3}{||c||}{<1,2,1>}\\
\multicolumn{3}{||c||}{<1,1,2>, <2,1,1>}\\
\hline\hline
\end{array}$$

$$\begin{array}{||c||}
\multicolumn{1}{c}{R_4} \\
\hline\hline
<1,1,2> = e \\
<2,1,1> = e \\
\hline\hline
\end{array}$$
$\pi({\bf A_4})$ is the free group generated by the single generator
$<1,2,1>$. This can be thought of as the fundemental group of the
noncommuting pentagon (=$\pi(S^1)$)

\vspace{2cm}

$$\begin{array}{||c|c|c|c||}
\multicolumn{4}{c}{G_5}\\
\hline\hline
\bf A_1 \times A_4 & \bf A_2 \times A_3 & \bf A_3 \times A_2 & \bf A_4
\times A_1
\\
\hline
\lam <1,2,1>   & & & \rho <1,2,1> \\
\hline
\multicolumn{4}{||c||}{<1,3,1>}\\
\multicolumn{4}{||c||}{<1,2,2>,<2,2,1>}\\
\multicolumn{4}{||c||}{<1,1,3>,<2,1,2>,<3,1,1>}\\
\hline\hline
\end{array}$$

$$\begin{array}{||c||}
\multicolumn{1}{c}{R_5} \\
\hline\hline
\\ \hline
<1,1,3> = e \\
<2,1,2> = e \\
<3,1,1> = e \\
\hline\hline
\end{array}$$
$\pi({\bf A_5})$ is the free group generated by five generators.
There are six pentagons in ${\bf A_5}$ but they are linked up
in such a way that there are only five generators.
This is similar to the fact that although the cube has six faces
(squares) there are only five generators i.e. only five faces must
commute in order for the entire cube to commute.

There is one unwritten trivial relation
$$<1,3,1> *  \rho \lam <1,1,1> * <1,3,1>^{-1} = \lam \rho <1,1,1>. $$
However since $X<1,1,1>$ is set to $e$, the relation is superfluous.
This relation corresponds to the center square in ${\bf A_5}$.
\vspace{2cm}

$$\begin{array}{||c|c|c|c|c||}
\multicolumn{5}{c}{G_6}\\
\hline\hline
\bf A_1 \times A_5 & \bf A_2 \times A_4 & \bf A_3 \times A_3 & \bf A_4
\times A_2 & \bf A_5 \times A_1
\\
\hline
\lam\lam <1,2,1>&\lam^2<1,2,1>& &\rho^2 <1,2,1>&\rho\lam <1,2,1>\\
\lam\rho <1,2,1> & & & & \rho\rho<1,2,1> \\
\lam <1,3,1> & & & & \rho <1,3,1> \\
\lam <1,2,2> & & & & \rho <1,2,2> \\
\lam <2,2,1> & & & & \rho <2,2,1> \\
\hline
\multicolumn{5}{||c||}{<1,4,1>}\\
\multicolumn{5}{||c||}{<1,3,2>,<2,3,1>}\\
\multicolumn{5}{||c||}{<1,2,3>,<2,2,2>,<3,2,1>}\\
\multicolumn{5}{||c||}{<1,1,4>,<2,1,3>,<3,1,2>,<4,1,1>}\\
\hline\hline
\end{array}$$

$$\begin{array}{||c||}
\multicolumn{1}{c}{R_6} \\
\hline\hline
<1,4,1> *  \rho \lam <1,2,1> * <1,4,1>^{-1} = \lam \rho <1,2,1> \\
<1,1,4> *  \lam^2 <1,2,1> * <1,1,4>^{-1} = \lam \lam <1,2,1> \\
<4,1,1> *  \rho\rho <1,2,1> * <4,1,1>^{-1} = \rho^2 <1,2,1> \\
\hline
<1,1,4> = e \\
<2,1,3> = e \\
<3,1,2> = e \\
<4,1,1> = e \\
\hline\hline
\end{array}$$
The first HNN-extension is the first non-trivial relation that we have.
The second relation is an amalgamation because $<1,1,4> = e$ and so
the relation reduces to
$$ \lam^2 <1,2,1>  = \lam \lam <1,2,1>.$$
This essentially shows the associativity of the free monoid on two
generators. Similarly for the third relation. So $\pi({\bf A_6})$ is
a group with 22 generators (12 old + 10 new.) Two pairs of
old
generators are set equal to each other and four new generators are
set
to $e$. We are left with 16 generators (10 old + 6 new.) There
is one nontrivial relation on these generators. However the group is
isomorphic to the free group on 15 generators.

\small
$$\begin{array}{||c|c|c|c|c|c||}
\multicolumn{6}{c}{G_7}\\
\hline\hline
\bf A_1 \times A_6 & \bf A_2 \times A_5 & \bf A_3 \times A_4 & \bf A_4
\times A_3 & \bf A_5 \times A_2 &\bf A_6 \times A_1
\\
\hline
\lam\lam^2 <1,2,1>&\lam^2\lam <1,2,1>&\lam^3<1,2,1> &\rho^3<1,2,1>
&\rho^2 \lam
<1,2,1>&\rho\lam^2 <1,2,1>
\\
\lam\lam\rho <1,2,1> &\lam^2 \rho <1,2,1>& & &\rho^2 \rho <1,2,1> &
\rho\lam\rho
<1,2,1>
\\
\lam\lam <1,3,1> &\lam^2 <1,3,1>& & &\rho^2 <1,3,1> & \rho\lam <1,3,1>
\\
\lam\lam <1,2,2> &\lam^2 <1,2,2>& & &\rho^2 <1,2,2> & \rho\lam<1,2,2>
\\
\lam\lam <2,2,1> &\lam^2 <2,2,1>& & &\rho^2 <2,2,1> & \rho\lam <2,2,1>
\\
\lam\rho^2 <1,2,1> & & & & & \rho\rho^2 <1,2,1>\\
\lam\rho\lam <1,2,1> & & & & & \rho\rho\lam <1,2,1>\\
\lam \rho <1,3,1> & & & & & \rho \rho <1,3,1>\\
\lam \rho <1,2,2> & & & & & \rho \rho <1,2,2>\\
\lam \rho <2,2,1> & & & & & \rho \rho <2,2,1>\\
\lam <1,4,1> & & & & & \rho <1,4,1>\\
\lam <1,3,2> & & & & & \rho <1,3,2>\\
\lam <2,3,1> & & & & & \rho <2,3,1>\\
\lam <1,2,3> & & & & & \rho <1,2,3>\\
\lam <2,2,2> & & & & & \rho <2,2,2>\\
\lam <3,2,1> & & & & & \rho <3,2,1>\\
\hline
\multicolumn{6}{||c||}{<1,5,1>}\\
\multicolumn{6}{||c||}{<1,4,2>,<2,4,1>}\\
\multicolumn{6}{||c||}{<1,3,3>,<2,3,2>,<3,3,1>}\\
\multicolumn{6}{||c||}{<1,2,4>,<2,2,3>,<3,2,2>,<4,2,1>}\\
\multicolumn{6}{||c||}{<1,1,5>,<2,1,4>,<3,1,3>,<4,1,2>,<5,1,1>}\\
\hline\hline
\end{array}$$

$$\begin{array}{||c||}
\multicolumn{1}{c}{R_7} \\
\hline\hline
<1,5,1> *  \rho \lam \lam <1,2,1> * <1,5,1>^{-1} = \lam \rho \lam
<1,2,1>
\\
<1,5,1> *  \rho \lam \rho <1,2,1> * <1,5,1>^{-1} = \lam \rho \rho
<1,2,1>
\\
<1,5,1> *  \rho \lam <1,3,1> * <1,5,1>^{-1} = \lam \rho <1,3,1> \\
<1,5,1> *  \rho \lam <1,2,2> * <1,5,1>^{-1} = \lam \rho <1,2,2> \\
<1,5,1> *  \rho \lam <2,2,1> * <1,5,1>^{-1} = \lam \rho <2,2,1> \\
\\
<1,4,2> *  \rho^2 \lam <1,2,1> * <1,4,2>^{-1} = \lam \rho^2 <1,2,1> \\
<2,4,1> *  \rho \lam^2 <1,2,1> * <2,4,1>^{-1} = \lam^2 \rho <1,2,1> \\
\\
***  <1,2,4> *  \lam^3 <1,2,1> * <1,2,4>^{-1} = \lam \lam^2 <1,2,1> \\
***  <4,2,1> *  \rho \rho^2 <1,2,1> * <4,2,1>^{-1} = \rho^3 <1,2,1> \\
<1,1,5>:\\
\lam^2 \lam <1,2,1> = \lam \lam^2 <1,2,1>\\
\lam^2 \rho <1,2,1> = \lam \lam \rho<1,2,1>\\
\lam^2 <1,3,1> = \lam \lam<1,3,1>\\
\lam^2 <1,2,2> = \lam \lam<1,2,2>\\
\lam^2 <2,2,1> = \lam \lam<2,2,1>\\
<5,1,1>:\\
\rho \rho \lam <1,2,1> = \rho^2 \lam <1,2,1> \\
\rho \rho^2 <1,2,1> = \rho^2 \rho <1,2,1> \\
\rho \rho <1,3,1> = \rho^2 <1,3,1> \\
\rho \rho <1,2,2> = \rho^2 <1,2,2> \\
\rho \rho <2,2,1> = \rho^2 <2,2,1> \\
<2,1,4>: \\
\lam^3 <1,2,1> = \lam^2 \lam <1,2,1>\\
<4,1,2>: \\
\rho^2 \rho <1,2,1> = \rho^3 <1,2,1> \\
\hline
<1,1,5> = e \\
<2,1,4> = e \\
<3,1,3> = e \\
<4,1,2> = e \\
<5,1,1> = e \\
\hline
\lam \# R_6: \\
\lam <1,4,1> * \lam \rho \lam <1,2,1> * \lam <1,4,1>^{-1} = \lam \lam
\rho <1,2,1> \\
\rho \# R_6: \\
\rho <1,4,1> * \rho \rho \lam <1,2,1> * \rho <1,4,1>^{-1} = \rho \lam
\rho <1,2,1> \\
\hline\hline
\end{array}$$
\normalsize

There are 57 (42 old + 15 new) generators. After
amalgamations
and setting some of the new generators to $e$ we have 42 (32
old + 10 new) generators. There are 11 relations on these 42
generators.

If you combine all the relations you get
$$ \rho^3 <1,2,1> = \rho \rho^2 <1,2,1> = \rho^2 \rho <1,2,1> $$ i.e.
$\rho$ is associative. Similarly for $\lam$.

This group has examples of the first nontrivial relations. Consider
$$<4,2,1> *  \rho \rho^2 <1,2,1> * <4,2,1>^{-1} = \rho^3 <1,2,1> $$
which reduces to
$$<4,2,1> *  \rho^3 <1,2,1> * <4,2,1>^{-1} = \rho^3 <1,2,1> $$
or
$$<4,2,1> *  \rho^3 <1,2,1> = \rho^3 <1,2,1> * <4,2,1>.$$
This relation can actually be seen. Look at the following
diagram of associations and reassociations.
\newpage
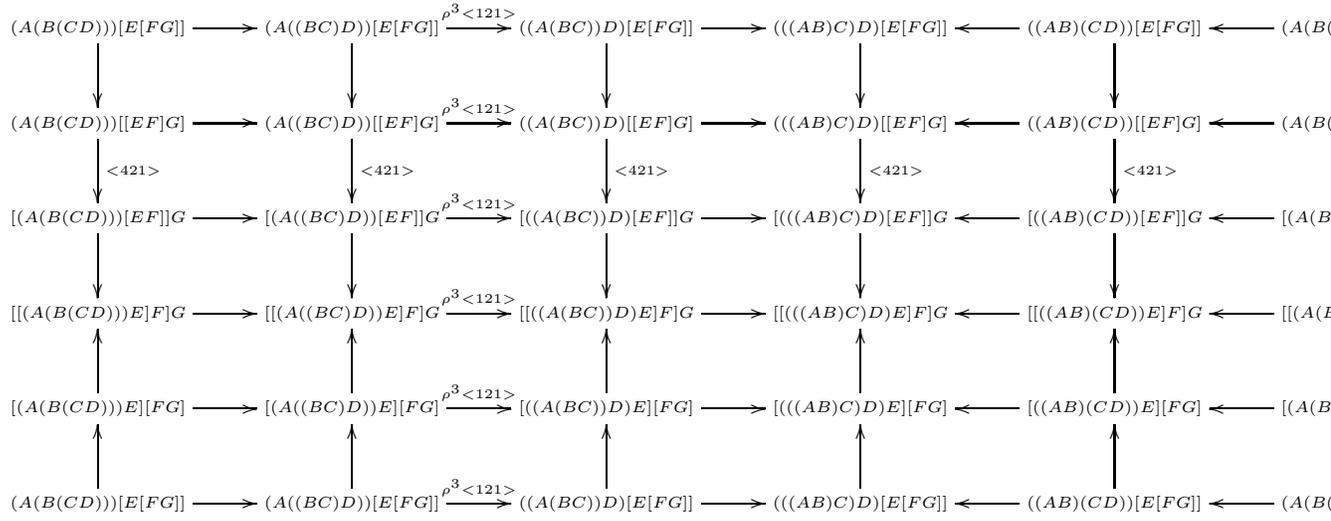
\tiny
$$ \begin{diagram}
(A(B(CD)))[E[FG]] \rto \dto& (A((BC)D))[E[FG]] \rto^{\rho^3<121>} \dto&
((A(BC))D)[E[FG]]
\rto
\dto &
(((AB)C)D)[E[FG]] \dto & ((AB)(CD))[E[FG]] \lto \dto & (A(B(CD)))[E[FG]]
\lto
\dto
\\
(A(B(CD)))[[EF]G] \rto \dto^{<421>}& (A((BC)D))[[EF]G]
\rto^{\rho^3<121>}
\dto^{<421>}&
((A(BC))D)[[EF]G]
\rto
\dto^{<421>} &
(((AB)C)D)[[EF]G] \dto^{<421>} & ((AB)(CD))[[EF]G] \lto \dto^{<421>}  & (A(B(CD)))[[EF]G]
\lto
\dto^{<421>}
\\
[(A(B(CD)))[EF]]G \rto \dto& [(A((BC)D))[EF]]G \rto^{\rho^3<121>}  \dto&
[((A(BC))D)[EF]]G
\rto
\dto &
[(((AB)C)D)[EF]]G \dto & [((AB)(CD))[EF]]G \lto \dto & [(A(B(CD)))[EF]]G
\lto
\dto
\\
[[(A(B(CD)))E]F]G \rto & [[(A((BC)D))E]F]G \rto^{\rho^3<121>}  &  [[((A(BC))D)E]F]G
\rto
&
[[(((AB)C)D)E]F]G &[[((AB)(CD))E]F]G\lto&[[(A(B(CD)))E]F]G\lto
\\
[(A(B(CD)))E][FG]\rto\uto&[(A((BC)D))E][FG]\rto^{\rho^3<121>}
\uto&[((A(BC))D)E][FG]\rto
\uto&
[(((AB)C)D)E][FG]\uto&[((AB)(CD))E][FG]\lto\uto&[(A(B(CD)))E][FG]\lto\uto
\\
(A(B(CD)))[E[FG]] \rto \uto& (A((BC)D))[E[FG]] \rto^{\rho^3<121>}  \uto&
((A(BC))D)[E[FG]]
\rto
\uto &
(((AB)C)D)[E[FG]] \uto & ((AB)(CD))[E[FG]] \lto \uto & (A(B(CD)))[E[FG]]
\lto
\uto
\end{diagram} $$
\normalsize

\newpage
The horizontal maps move the round parenthesis and the vertical maps
move the square parenthesis
The top row and the bottom row are the same. As are the left side and
the right
side. Each square commutes out of naturality. The whole
thing is a torus. As with every torus, the generators commute. The two
generators have been marked off. Going around $\rho^3 <1,2,1>$ and then
going around $<4,2,1>$ is the same as doing it vice versa. All other
maps are part of equivalence classes that have a $1$ in the center
position of the $<i,j,k>$ bracket and hence are set equal to $e$.

The important point is that the two stared relations show commutativity
of generators. These generators do not show up in any other relations
and hence the group is not free!
The higher groups contain this group and hence they are also not free.

\section{Fundamental groups of quotients}
In order to look at quotients of associative categories, we must define
congruences that respect the associative category structure.
\begin{defi}
A congruence on an associative category $(\bf B, \otimes, \beta)$ is an
equivalence relation $\sim$ for each Hom set of morphisms of $\bf B$
such that
\begin{enumerate}
\item composition is respected;
\item the $\otimes$ is respected, i.e. if $f \sim f'$ and $g \sim
g'$, then $f \otimes g \sim f' \otimes g'$;
\item naturality is respected, i.e. for any $f:A \longrightarrow A'$, $
 g:B \longrightarrow B'$ and $h:C \longrightarrow C'$, we have
$$ ((f \otimes g) \otimes h) \circ \beta_{A,B,C} \quad \sim \quad
\beta_{A',B',C'} \circ (f \otimes (g \otimes h)).$$
\end{enumerate}
\end{defi}
By the customary universal algebra tricks, to every congruence on
$(\bf B, \otimes, \beta)$ there exists a unique associative category
$({\bf \tilde{B}},
\tilde{\otimes}, \tilde{\beta})$ and a unique strict tensor functor
$\Pi:\bf B
\longrightarrow {\bf \tilde{B}}$ that satisfies universal properties.
${\bf \tilde{B}}$ is isomorphic on objects. Morphisms are equivalence
classes of morphisms of $\bf B$, denoted $[f]$. The tensor product is
defined as usual:  $[f] \tilde{\otimes} [g] = [f \otimes g]$. The
associativity is $\tilde{\beta}_{A,B,C}=[\beta_{A,B,C}]$. Condition
3 insures naturality of the associativity.

The next step is to look at quotients of the 2-sketch $\bf A$ of the
theory of associative categories. Since $\bf A$ is a strict tensor
2-category, we must define what congruences of such categories are
like. Our only
interest is when 2-cells of $\bf A$ (which correspond to 1-cells of
$\Abar$) are set equal to
each other, not 0-cells or 1-cells. So we have the following definition.
\begin{defi}
A 2-congruence on a strict tensor 2-category $(\bf C, \otimes)$
is an equivalence relation $\sim$ on each set of 2-cells between any
two 1-cells of $\bf C$ such that
\begin{enumerate}
\item the $\otimes$ is respected, i.e. if $\phi \sim \phi'$ and $\gamma
\sim \gamma'$, then $(\phi \otimes \gamma) \sim (\phi' \otimes
\gamma')$;
\item vertical composition is respected, i.e. if $\phi \sim \phi'$ and
$\gamma
\sim \gamma'$, then $$(\phi \circ_V \gamma) \sim (\phi' \circ_V
\gamma');$$
\item horizontal composition is respected (similar to 2).
\end{enumerate}
\end{defi}

By generalized universal algebra, to every 2-congruence on a strict
tensor
2-category
$(\bf C, \otimes)$ there exists a unique strict tensor 2-category
$({\bf
\tilde{C}},
\tilde{\otimes})$ and a unique strict tensor 2-functor $\Pi:\bf C
\longrightarrow {\bf \tilde{C}}$ that satisfies the usual universal
properties.
${\bf \tilde{C}}$ is isomorphic on 0-cells and 1-cells. 2-cells
are equivalence
classes of 2-cells of $\bf C$, denoted $[\phi]$. The tensor
product, vertical composition and horizontal composition
are defined as usual. The ``commutativity'' of vertical composition and
horizontal composition in
${\bf \tilde{C}}$
, falls out of the fact that $\bf C$ has this
property. The converse of this statement is also true: given two 
strict tensor 2-categories and a strict tensor
2-functor between them that is isomorphic on 0-cells and 
1-cells, we get such a
2-congruence.

We shall now examine what a 2-congruence for our 2-sketch
$\bf
A$ is like. The fact that the 2-congruence preserves tensor product
means
that if $\phi=(\phi_1, \ldots , \phi_k)  \sim \phi'=(\phi'_1, \ldots ,
\phi'_{k'})$ and
$\gamma=(\gamma_1, \ldots , \gamma_l)  \sim (\gamma'_1,
\ldots,\gamma'_{l'}) = \gamma'$ then
$$\phi \otimes \gamma = (\phi_1, \ldots , \phi_k, \gamma_1, \ldots ,
\gamma_l)  \sim (\phi'_1, \ldots ,\phi'_{k'},\gamma'_1, \ldots ,
\gamma'_{l'}) = \phi' \otimes \gamma'.$$
Since we are interested in
associative categories in ${\bf Cat}$ i.e. strict tensor 2-functors $F$
from $(\bf A, \otimes)$ to $({\bf Cat},\times)$, we have the following
$$F(\phi) \times F(\gamma) =
F(\phi \otimes \gamma) =
F(\phi' \otimes \gamma') =F(\phi') \times F(\gamma'), $$ and hence
using usual properties of the product in ${\bf Cat}$
$F(\phi)= F(\phi')$ and $F(\gamma) = F(\gamma')$. Such functors
determine quotients and hence if $(\phi \otimes \gamma) \sim (\phi'
\otimes \gamma')$ then $\phi \sim \phi'$ and $\gamma \sim \gamma'.$
By a small inductive argument it suffices to discuss only individual
components of the 2-cells. Hence we shall talk about $\phi_i \sim
\phi'_i$ where they are both elements of ${\bf A_{n_i}}$ for some $n_i$.
Hence we have, that every 2-congruence on $\bf A$ induces a congruence
on each of the $\An$.

Vertical composition of 2-cells in $\bf A$ is simply composition in each
${\bf A_{n_i}}$.

Finally we are left with horizontal composition which is the heart of
the matter. We remind the reader that horizontal composition is defined
with the use of the operad $Q$. A review of how the operad is defined
and employed would be beneficial at this time. If
$\phi \sim \phi'$ and $\gamma \sim \gamma'$ then
$$\Lambda = (\phi \circ_H \gamma) \sim (\phi' \circ_H \gamma') =
\Delta.$$ By the argument above, this means that
$$\Lambda_i = Q(\phi_i,\gamma_{n_1},\ldots , \gamma_{n_s}) \sim
Q(\phi'_i,\gamma'_{n_1},\ldots , \gamma'_{n_s}) = \Delta_i $$ where the
subscripts of $Q$ are abandoned since they are the same $Q$.

Between every two 1-cells of $\bf A$ there is a distinct 2-cell, denoted
$\breve{\gamma}$, that has all its components $\breve{\gamma_i}$ in the
MIS ${\bf T_{n_i}}$. We call such a 2-cells a``branch'' (branches make
up trees.) Given two branches $\breve{\gamma}$ and
$\breve{\gamma'}$, their tensor product $\breve{\gamma}
\otimes \breve{\gamma'}$ is also a branch. Similarly for vertical
composition. ``Branchness'' is also preserved by horizontal composition
because it is preserved by the operad $Q$ i.e.
$$Q(\breve{\phi},\breve{\gamma_{n_1}}, \ldots ,\breve{\gamma_{n_s}}) =
\breve{\phi'}.$$ This is shown by an inductive proof on the construction
of the $Q$'s and from the fact that the $V_{i,j,k}$ are built up from
earlier $V_{i,j,k}$.

For every congruence, there is a unique equivalence class
between every two 1-cells, namely, those 2-cells congruent to the unique
branch 2-cell. We denote this equivalence class by $[\breve{\phi}]$.

Given an associative category $(\bf B, \otimes, \beta)$, there are
unique strict tensor functors $P_n:\An \longrightarrow \Bitn$ as was
shown in Section 4. We use $P_n$ in the following diagram.
$$ \begin{diagram}
\Tn \rrto^{!} \ddto_{L_n}
& & \bf 1 \ddto \ar@{=}[rr]& & \bf 1 \ar[dddd]
\\ \\
\An \rrto \ar[dd]_{P_n}
& &
\pi(\An) \ar@{-->}[ddrr]^{\pi(P_n)}
\\ \\
\Bitn \ar[rrrr] & & & & \pi(\Bitn).
\end{diagram} $$
Both squares are pushouts.
The inner square is familiar from Section 7. The result of the outer
pushout is $\pi(\Bitn)$ which is shortened to $\pi(\bf B_n)$
(this should not seem so strange since the tensor product and the
reassociation of $\Bitn$ is basically the same as the tensor product and
the reassociation of $\bf B$.) The induced arrow is denoted by
$\pi(P_n)$. The functorial properties of $\pi(-)$ will not be discussed
here. Since the left hand vertical map is surjective on objects,
$\pi({\bf B_n})$ is a group and we call it the fundamental group of $\bf
B_n$. We must stress that $P_n$ is not necessarily full (surjective on
morphisms). A typical example is when $\bf B$ is a {\em strict} tensor
category with $\beta$ being something other then the identity. Then $\bf
A_3$ is the indiscrete category with two objects whereas the ${\bf
It(B)_3}$ has one object and - by composition - infinite morphisms. Our
interest lies not in $\pi(\bf B_n)$ but in the image of $\pi(P_n)$ since
there are to be found the morphisms that are of concern to
coherence.

Now for a short discussion of the way $Q$ works. $Q(\phi,\gamma, \ldots,
\gamma)$ takes $\phi$ and for each letter that $\phi$ reassociates,
puts in other letters. How many other letters? One or more, depending on
what partition we have. Now if $\phi$ is a generating morphism of the
form $\phi= V_{i,j,k}(\phi_i, \phi_2, \phi_3)$ then
$$Q(\phi,\gamma, \ldots,
\gamma) = V_{i',j',k'}(\phi'_i, \phi'_2, \phi'_3)$$ by the definition of
$Q$ given in Section 1.4. This can be denoted in our group theoretic
notation as
$$Q(<i,j,k>, \gamma, \ldots, \gamma) = <i',j',k'>$$ where $\i' \ge i,
j'\ge j$ and $k' \ge k$. Similarly if $\phi = \phi_a \otimes \phi_b$
then
$$Q(\phi,\gamma, \ldots,
\gamma) = Q(\phi_a, \gamma, \ldots \gamma) \otimes
Q(\phi_b, \gamma, \ldots \gamma).$$
In group theoretic notation we have
$$Q(X<i,j,k>, \gamma, \ldots, \gamma) =X'<i',j',k'>$$ where $\i' \ge i,
j'\ge j, k' \ge k$ and $length(X') \ge length(X)$.

The important point is that every 2-cell (of $\bf A$ or 1-cell of $\An$)
that is a branch, goes to $e$, the identity of $\pi(\bf B_n)$.

Putting this all together we have the following.

\begin{teo} If $j \ge 2$, $i + j + k = n$ and $\pi(P_n)(<i,j,k>) = e$ (
the
identity of the $\Bitn$ ), then for all $i' \ge i$, $j' \ge j$ and $k'
\ge k$ with $i' + j' + k' = n'$ we have $\pi(P_{n'})(<i',j',k'>)=e$.
\end{teo}

This theorem is a proper generalization of the Mac Lane's coherence
theorem and we have:

\begin{cor}[Mac Lane's Coherence Theorem] If $\pi(P_4)(<1,2,1>)=e$,
then for all $i'$, $j'$ and $k'$ with $i' + j' + k' = n'$ we have
$\pi(P_{n'})(<i',j',k'>)=e$ and hence $\pi(P_{n'})$ is the trivial map
for all
$n'$.
\end{cor}
{\bf Proof.} Set $i=1, j=2$ and $k=1$ in the above theorem. The fact
that $\pi(P_{n'})$ is trivial for all $n'$ comes from the fact that all
the generators of $\pi(\An)$ are of the form $X<i,j,k>$ where $X$ is an
element of the free monoid on $\lam, \rho$. If $<i,j,k>=e$, then
$X<i,j,k> = Xe = e \Box$

Let us return to a typical category that was briefly discussed in 
Section 2. Let $R$ be a ring with a unit and ${\bf R-Mod}$ 
denote the category of left and right $R$ modules. ${\bf R-Mod}$ has 
the usual tensor product. Then we have the set of all nonpathalogical 
reassociation natural
transformations $$\beta_{A,B,C} : A \otimes
(B \otimes C) \longrightarrow (A \otimes B)
\otimes C $$ of the form $$ a \otimes( b \otimes c) \longmapsto \phi(a' \otimes b') \otimes c'.) $$
Since $\beta$ must be natural, $\phi$ is a central element of the ring.
Since $\beta$ must be an isomorphism, $\phi$ must have a multiplicative
inverse. So $\phi \in CU(R)$, the central units of $R$. $CU(R)$ forms an 
abelian multiplicative group. Let $\phi \in CU(R)$, then $O(\phi) \in \bf N \cup\{ \infty\}$
is the number of times that you must go around the pentagon in order to  get
to the identity. If $O(\phi)=1$ then $\phi$ induces a coherent reassociator.
This order turns out to be nothing more then asking what power do
we have to raise $\phi$ to, in order to get the unit
i.e. $O(\phi)=n$ iff $\phi^n=1$. So noncoherence for ${\bf It(R-Mod)_4}$ is
classified by $CU(R)$.

The conclusion is that there are many forms of noncoherence but we can
classify them as a tree. If going around a certain loop in
$\An$ is set equal the identity, then there are implications for the
higher $\An$. That is, every relation in $\An$ has ramifications 
for the ${\bf A}_k$ for $k \ge n$. However, for non-unital
associative categories, there are no ramifications for ${\bf A}_l$ for
 $l \le n$. When we add the unit, we make the classification more interesting.

\section{Quotients of unital associative categories}
Since we have not given a formal construction of $\Abar'$ (see Section
6,) the free
associative category with a unit or of ${\bf A'}$ the theory of
associative categories with units, we shall give an intuitive rather
than a formal discussion
of their fundamental groups. All congruences must also take into 
account the morphisms
$L,R :\An \longrightarrow {\bf A_{n-1}}$. Hence there are, in fact,
fewer congruences if we insist upon units. Here is another way of
looking at it.  Assume there is a loop in
$\Bitn$: we can write it as
$$\xymatrix{
f_1 \ar[r]^{\phi_1} & f_2 \ar[r]^{\phi_2} & \ldots \ar[r] & f_n
\ar[r]^{\phi_n} & f_1} $$
where the $f_i$'s are associations of any $n$ objects in $\bf B$. In
particular, one of the objects may be the unit object $I$. Let $f_i^*$
be
the association $f_i$ without that unit object. $f_i^*$ is then an
association of $n-1$ objects. Hence we have the following commutative
diagram:
$$\xymatrix{
\Bitn & f_1 \ar[r]^{\phi_1} \ar[dd] & f_2 \ar[r]^{\phi_2} \ar[dd] &
\ldots
\ar[r] \ar[dd]
& f_n
\ar[r]^{\phi_n} \ar[dd] & f_1 \ar[dd] \\ \\
{\bf It(B)_{n-1}}
& f_1^* \ar[r]^{\phi_1^*} & f_2^* \ar[r]^{\phi_2^*} & \ldots \ar[r] &
f_n^*
\ar[r]^{\phi_n^*} & f_1^*.} $$ where the vertical maps are identities
tensored with $L$ or $R$. The squares commute only if we assume the
identity triangle coherence requirement discussed in Section
6. If the top loop is the identity, then the bottom loop must also be
the identity. In other words, when there are no units a
relation in $\An$ causes a relation in ${\bf A_{n+t}} $. In contrast,
when there are units, relations in ${\bf A_{n+t}}$ {\em go down} to
the lower $\An$. Hence there are fewer congruences possible.
So, noncoherence for unital associative categories are classified 
by $\bf Z$ the integers. To every unital associative category, there 
is an integer $v$ assigned to it. $v$ is the least number that one must 
go around a generator of some -- hence all -- $\An$ to get the identity.

\section{Future Directions}

There are several places that the work done here might be applied.
Most of the second half of \cite{Shn&Stern} is dedicated to calculating
the cohomology of Drinfeld algebras (bialgebras that are coassociative
only up to a coherent isomorphism.) It would be interesting to study the
relationship of these cohomology groups to the cohomology of bialgebras
that are coassociative only up to an isomorphism - without a coherence
requirement. Such bialgebras arise when representing (in the sense of,
say, \cite{Yetter} )  an associative category. We also feel that the
groups presented in section 8 will be of importance to the
study
of pentagonal algebras and homotopy associative (Lie) algebras (see
\cite{Stasheff} .)

The second part of \cite{Thesis} deals with functors between associative
categories that do not neccesarily satisfy the hexagon coherence
condition (see \cite{Epstein}.) An induction scheme to calculate the
fundamental groups of ``mapping funnals'' is given. A
generalisation of Proposition 4 in Section 5 is proved where the
restriction to strict tensor functors is lifted. There is an attempt to
classify noncoherent tensor functors. Much work remains to be done in
this area. The standard coherence theorem states that any coherent
tensor category category is equivalent to a strict tensor category via a
coherent tensor functor. Given a noncoherent tensor category can we
say that it is equivelent to a strict tensor category via a (weaker)
noncoherent tensor functor?

Many
of the assumptions made in this work can be relaxed for more interesting
computations. For example, in this paper, the reassociations are
always considered to be {\it iso}morphisms. What happens if we relax
this
requirement and only ask for a morphism? We would then be working with
general categories rather than (Catalan) groupoids. Would we get a
fundamental group or a fundamental {\it monoid}? Similarly for
noncoherent tensor functors (\cite{Thesis}), we assumed there is an
isomorphism
$F(A)
\otimes ' F(B) \longrightarrow F(A \otimes B)$. What happens if we
loosen this requirement and ask only for a morphism between them?
These situations are known to arise ``in nature''. For instance,
\cite{Yetter} is concerned with braidings that are not necessarily
isomorphisms (termed ``pre-braidings''.) Also, when using the forgetful
functor, $U: {\bf R-Mod} \longrightarrow {\bf Ab}$ from the module
category of an arbitrary commutative ring $R$ to the category of abelian
groups, the morphism $U(A) \otimes U(B) \longrightarrow U(A \otimes B)$
is in general not an isomorphism.

The next coherence requirement that is under investigation
is commutativity \cite{Yano}. This area is of utmost importance to the
study
of quantum groups, quasi-triangular Hopf algebras, quantum field theory,
etc.. \cite{Joy&Str} has given a three-level hierarchy of coherence
requirements for commutativity: symmetric, balanced and braided.
Each of these coherence conditions correspond to different algebraic
structures. Is this hierarchy complete? Are there intermediate levels?
Each one of these coherence rules corresponds to ``filling in'' part of
the permutoassociahedrons (\cite{Kapranov} ). We would like to look
at the fundamental groups of each of the related polytopes and see if
there are any other interesting coherence conditions.

Another area of extreme interest is categorical logic and coherence. We
would like to look more carefully at the coherence requirements for
cartesian closed categories. It has been shown (e.g. \cite{Mac Lane} for
a general survey, \cite{Kelley1}) that these coherence
requirements are intimately related to the cut-elimination theorem which
is central to proof theory of (intuitionistic) logic. The goal would be
to furnish a classification of categories that fail this coherence
condition and hence to see what can be learned about logical systems in
which the cut-elimination theorem fails. Coherence requirements have
also shown up in the area of linear logic. Linear logic
deals with -- the more general -- {\it monoidal}  closed categories (see
e.g. \cite{Troelstra} ).

There are numerous other coherence problems which we can explore using
generalizations of methods used in this paper. For example, we can
look at Laplaza's distributivity categories, Shum's tortile tensor
categories, Crane and Yetter's Hopf categories etc..

The long term goals are to study n-categories as forms of
higher dimensional algebras. Recently, these algebras have become
quite popular
Even mathematical
physicists (see e.g. \cite{Cra&Yet} or \cite{Boez&Dolan} )have
shown an interest
in such structures. In 1972, Kelley (\cite{Kelley2})
formulated the notion of a club in order to present coherence problems.
A club is like a presentation of a universal algebraic theory. However,
instead of just having operations and identities (or generators and
relations), there are {\it two} levels of operations and identities.
Loosely speaking, coherence requirements are second dimensional
identities. (From this point of view,  the Catalan groupoids constructed
in here are the 2-dimensional versions of the Cayley graphs of
groups.) In the thesis, we study categories with structures that are
``free'' or partially ``free'' of coherence requirements. In regular
(1-dimensional) algebras, the notions of a free algebra plays a major
role in homological algebra. Perhaps noncoherent categories will play
such a role in higher dimensional
homological algebra.

\end{document}